\documentclass[3p,times,sort&compress]{elsarticle}
\graphicspath{{./figures/.}}

\usepackage{amssymb}
\usepackage{physics}
\usepackage{subcaption}

\newcommand{\autocite}{\cite}
\newcommand{\textcite}{\cite}

\newcommand{\dxi}{\dd{\xi}}
\newcommand{\ds}{\dd{s}}
\newcommand{\dx}{\dd{x}}

\newcommand{\OO}{\mathcal{O}}
\newcommand{\ppi}{\mathrm{\pi}}
\newcommand{\e}{\mathrm{e}}
\newcommand{\dr}{\dd{r}}

\newcommand{\sfrac}[2]{{\textstyle\frac{#1}{#2}}}

\journal{Physica D}
\usepackage {tikz}
\usetikzlibrary {positioning}
\usetikzlibrary {calc}
\newcommand\getcurrentref[1]{%
 \ifnumequal{\value{#1}}{0}
  {??}
  {\the\value{#1}}%
}

\def\XXint#1#2#3{{\setbox0=\hbox{$#1{#2#3}{\int}$}
     \vcenter{\hbox{$#2#3$}}\kern-.5\wd0}}

\usepackage{mathtools}

\DeclareMathOperator\erfc{erfc}
\usepackage{hyperref}
\renewcommand{\i}{\mathrm{i}}
\begin{document}

\begin{frontmatter}

\title{Burgers' equation in the complex plane}

\author[qut]{Daniel J. VandenHeuvel}
\author[mac]{Christopher J. Lustri}
\author[not]{John R. King}
\author[qut]{Ian W. Turner}
\author[qut]{Scott W. McCue \corref{cor1}}
\address[qut]{School of Mathematical Sciences, Queensland University of Technology, Brisbane QLD 4001, Australia}
\address[mac]{School of Mathematical and Physical Sciences, 12 Wally's Walk, Macquarie University, New
South Wales 2109, Australia}
\address[not]{School of Mathematical Sciences, University of Nottingham, Nottingham NG7 2RD, United Kingdom}
\cortext[cor1]{Corresponding author: scott.mccue@qut.edu.au}


\begin{abstract}
Burgers' equation is a well-studied model in applied mathematics with connections to the Navier-Stokes equations in one spatial direction and traffic flow, for example.  Following on from previous work, we analyse solutions to Burgers' equation in the complex plane, concentrating on the dynamics of the complex singularities and their relationship to the solution on the real line.  For an initial condition with a simple pole in each of the upper- and lower-half planes, we apply formal asymptotics in the small- and large-time limits in order to characterise the initial and later motion of the singularities.  The small-time limit highlights how infinitely many singularities are born at $t=0$ and how they orientate themselves to lie increasingly close to anti-Stokes lines in the far field of the inner problem.  This inner problem also reveals whether or not the closest singularity to the real axis moves toward the axis or away.  For intermediate times, we use the exact solution, apply method of steepest descents, and implement the AAA approximation to track the complex singularities.  Connections are made between the motion of the closest singularity to the real axis and the steepness of the solution on the real line.  While Burgers' equation is integrable (and has an exact solution), we deliberately apply a mix of techniques in our analysis in an attempt to develop methodology that can be applied to other nonlinear partial differential equations that do not.
\end{abstract}




\end{frontmatter}



\section{Introduction}

Studying solutions of nonlinear partial differential equations (pdes) in the complex plane is a fascinating topic. Here we are interested in Burgers' equation,\bgroup\everymath{\displaystyle}
\begin{align}\label{eq:maineq}
\begin{cases}
\begin{aligned}
\pdv{u}{t} + u\pdv{u}{x}&= \mu \pdv[2]{u}{x}, \quad x \in \mathbb R,\,t>0,\\[5pt]
u(x, 0)&= u_0(x),\quad x\in\mathbb R,
\end{aligned}
\end{cases}
\end{align}\egroup
where $\mu > 0$. It is known in the case $\mu = 0$ that solutions to Burgers' equation can become increasingly steeper and evolve to form shocks in finite time, caused by (branch-point) singularities of the solution in the complex plane $x \in \mathbb C$, evolving towards and then ultimately touching the real axis. These observations were explored by Bessis \& Fournier \textcite{bessis1984pole, bessis1990complex}, for example, who analysed the solution in the complex plane as a way of understanding the blow-up of the (slope of the) solution on the real line. For $\mu > 0$, there is no such blow-up, since the steepening of the solution is overwhelmed by diffusion.  In this regime, insight into the behaviour of solutions in the complex plane has been progressed using a combination of analytical and numerical means \autocite{senouf1997dynamics, senouf1997dynamicsA,bessis1984pole, bessis1990complex, caflisch2015complex, weideman2021dynamics, Chapman_2007, weideman2003computing}. These ideas will be continued in this work and expanded upon, with our primary focus being the evolution of singularities in the solution of \eqref{eq:maineq} in the complex plane, and how these singularities affect the solution on the real line.

In terms of applications, Burgers' equation is a very well studied prototype model that demonstrates a competition between nonlinear wave steepening due to advection and linear smoothing due to diffusion.  It can be thought of as a crude simplification of the Navier-Stokes equations in one spatial dimension in the absence of a pressure gradient or body force.  In this context, $u$ is fluid velocity and $\mu$ is the fluid viscosity (we shall continue to refer to $\mu$ as the viscosity and the $\mu=0$ equation as the inviscid model).  Alternatively, Burgers' equation is considered as a simple model for traffic flow, again in one spatial direction.  In this case, if $\rho$ is the density of the traffic, then a reasonable constitutive relationship is that the flux is $\rho(1-\rho)-\mu \rho_x$, where here $\mu$ is a measure of how traffic slows due to a gradient in traffic ahead.  Conservation of mass, together with a change of variable $u=1-2\rho$, leads to Burgers' equation.

An obvious advantage of studying Burgers' equation \eqref{eq:maineq} is that it has an exact solution,
\begin{equation}\label{eq:exactsolution} u(x, t) = \dfrac{\displaystyle\int_{-\infty}^\infty \dfrac{x-s}{t}\exp\left\{-\dfrac{1}{2\mu}\left[\displaystyle\int_0^s u_0(\xi) \dxi + \dfrac{(x-s)^2}{2t}\right]\right\} \ds}{\displaystyle\int_{-\infty}^\infty \exp\left\{-\dfrac{1}{2\mu}\left[\displaystyle\int_0^s u_0(\xi) \dxi + \dfrac{(x-s)^2}{2t}\right]\right\} \ds},
\end{equation}
derived by relating \eqref{eq:maineq} to the heat equation through the Cole-Hopf transformation \cite{cole1951,hopf1950}.  Therefore, strictly speaking, we can compute the solutions at any $x$ and $t$, including complex values of $x$. However, numerical evaluation of the integrals in \eqref{eq:exactsolution} poses significant challenges in the complex plane, for example near singularities, or in various parameter combinations for which the integrands are highly oscillatory, including small $t$ and $\mu$ and values of $x$ far from the real axis. Therefore we shall also resort to asymptotic and other techniques to explore the evolution in the complex plane. Furthermore, while Burgers' equation \eqref{eq:maineq} may have an exact solution, most other nonlinear pdes do not, and we are therefore motivated to employ techniques that can generalise to other nonlinear models.

It is easy to see that singularities of solutions to \eqref{eq:maineq} in the complex plane must be simple poles. Suppose such a singularity occurs at $z = s(t)$ (where we replace the variable $x$ with $z$ to emphasise that it is complex), with
\begin{equation}\label{eq:poleansatz}
u(z, t) \sim \frac{a(t)}{\left(z-s(t)\right)^b}\quad\text{as}\quad z \to s(t),
\end{equation}
then, provided $b$ is not $0$ or $-1$, to leading order the dominant balance in \eqref{eq:maineq} comes from $uu_x$ and $\mu u_{xx}$, which implies $b=1$ and $a(t)=-2\mu$. That is, all singularities are simple poles with residues $-2\mu$ (this result is a step towards the claim that Burgers' equation possesses the Painlev\'{e} property for pdes~\cite{weiss1983}).
By including a correction term, this argument  can be extended to show
\begin{equation}
u(z, t) \sim -\frac{2\mu}{z-s(t)}+ \dv{s}{t}\quad\text{as}\quad z \to s(t).
\label{eq:localpole}
\end{equation}
For the other relevant case, namely $b=-1$ in \eqref{eq:poleansatz}, the leading-order balance involves only $u_t$ and $\mu u_{xx}$, corresponding to simple zeros.

In their early study, Bessis \& Fournier~\textcite{bessis1984pole, bessis1990complex} considered solutions of \eqref{eq:maineq} using the cubic initial condition $u(x, 0) = 4x^3 - x/t_s$, where $t_s > 0$ is the blow-up time for the corresponding solution when $\mu = 0$. Their work for $\mu>0$ was later extended in some detail by Senouf \textcite{senouf1997dynamics, senouf1997dynamicsA} using the same initial condition, giving much deeper understanding of the dynamics of the poles, especially in the limiting case $\mu\ll 1$.  Following early work by Sulem et al.~\cite{sulem1983}, both Caflisch et al.~\textcite{caflisch2015complex} and Weideman \textcite{weideman2021dynamics} consider the qualitatively similar initial condition $u(x, 0) = -\sin(x)$ and study the resulting complex singularities using a variety of analytical and numerical approaches. While all of the above studies are certainly complicated in many respects, these two initial conditions lead to a slightly less challenging analysis in the sense that the resulting singularities are confined to the imaginary axis for all time (due to the initial data -- and hence the solution for all $t$ --  being odd in $x$), simplifying some of the tools required for progress to be made. Moreover, these initial conditions are both entire functions, meaning that all poles come in from infinity at $t=0$ and evolve towards the real axis (before turning around and moving away).  As such, these studies avoid dealing with singularities emerging from within the plane at $t=0$.

These considerations motivate us to consider a specific initial condition whose analytic continuation is not an entire function and can serve as a prototype for more general initial conditions. One reason for this interest in initial conditions which are not entire is that infinitely many poles will emerge from the singularities of the initial condition, and this situation has attracted much less interest in the literature. An  obvious choice is an initial condition with a single simple pole in each of the lower and upper halves of the complex plane, giving
\begin{equation}\label{eq:ic}
u(x, 0) = u_0(x) = \frac{1}{1+x^2},
\end{equation}
which has simple poles at $x = \pm \i$. Our work complements that of Chapman et al.~\textcite{Chapman_2007} who consider the same initial condition \eqref{eq:ic} and apply exponential asymptotic methods to the study of the resulting solution $u(x, t)$ in the complex plane in the limit $\mu\rightarrow 0^+$, both with and without reference to the exact solution \eqref{eq:exactsolution}. Our work here differs from these authors in a number of ways. While Chapman et al.~\textcite{Chapman_2007} study (\ref{eq:maineq}) with \eqref{eq:ic} in considerable detail, their primary focus was on the limit $\mu \to 0^+$ and the associated near-shock behaviour, which we avoid.  In contrast, our concern is for a wide range of values of viscosity $\mu$, in particular the case $\mu = \OO(1)$.  Having said that, we shall present some results for the case $\mu\ll 1$ that complement those for $\mu = \OO(1)$. Note that the emergence of infinitely many singularities from one initial singularity is a phenomenon of more general interest in the study of nonlinear pdes and moving-boundary problems, and is, for example, reminiscent of problems associated with Hele-Shaw flow \cite{tanveer1993,tanveer2000,costintanveer2006}.

The inviscid version of (\ref{eq:maineq}) (i.e., with $\mu=0$), together with \eqref{eq:ic}, is straightforward to solve using the method of characteristics.  On the real line, the solution profile steepens until a shock forms at $(x_s,t_s) = (\sqrt{3}, 8\sqrt{3}/9)$.  In the complex plane, the solution has branch point singularities where $z=z_s(t)$ and $u=u_s^*(t)$ are solutions of
\begin{equation}\label{eq:blowupconditions}
u_s^* = \frac{1}{1+\left(z_s-u_s^*t\right)^2}, \quad 1 = \frac{2t\left(z_s-u_s^*t\right)}{\left[1+\left(z_s-u_s^*t\right)^2\right]^2}.
\end{equation}
Analysis of these equations leads to a quartic polynomial and indicates that the solution has two branch points in each of the upper and lower half planes.  The two branch points that are closest to the real axis move towards this axis and then touch it at precisely the time that the shock forms on the real line, namely $(x_s,t_s) = (\sqrt{3}, 8\sqrt{3}/9)$.  More details of this analysis are given in the Supplementary Material (and in~\textcite{Chapman_2007}), where we also show that, for small time, the two branch points in the upper-half plane leave $z=\i$ as
\begin{equation}\label{eq:inviscidsmall}
z_s \sim \i \pm (1-\i)t^{1/2} + \sfrac{1}{4}(2 \pm 1)t, \quad t \to 0^+.
\end{equation}

\begin{figure}[!h]
\centering
\includegraphics[width=\textwidth]{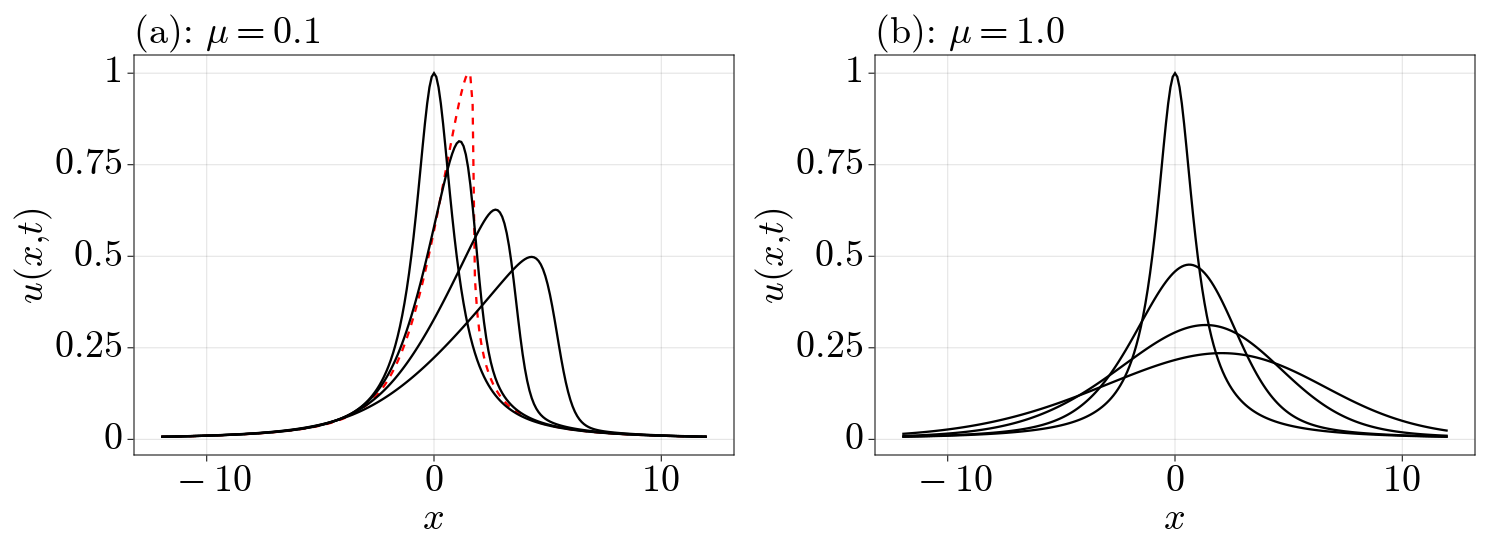}
\caption{Profiles of solutions of Burgers' equation (\ref{eq:maineq}) with (\ref{eq:ic}) for increasing time (black lines), computed for (a) $\mu=0.1$ and (b) $\mu=1$.  In each image, the profiles are shown for times $t = 0$, $t_s \approx 1.5396$, $5$, and $10$. In (a), the (red) dashed curve is the solution for $\mu=0$ at $t=t_s$.
}
\label{fig:compareinviscidviscous}
\end{figure}

It is well known that for $\mu>0$, the shock is regularised so that the viscous solution exists for all time $t>0$.  This regularisation is illustrated graphically in figure~\ref{fig:compareinviscidviscous}. Figure~\ref{fig:compareinviscidviscous}(a) is for $\mu = 0.1$, which is a representative small value of $\mu$.  The solution of (\ref{eq:maineq}) with \eqref{eq:ic} is provided by (black) solid curves for four times $t=0$, $8\sqrt{3}/9$, 5, and 10.  While the profile begins to steepen as time increases from $t=0$, there comes a time at which it begins to flatten.  Clearly the solution continues to exist for $t>8\sqrt{3}/9$, which is the time that the shock forms for $\mu=0$. For comparison, the solution with $\mu=0$ at $t=8\sqrt{3}/9$ is also included as a (red) dashed curve. Figure~\ref{fig:compareinviscidviscous}(b) is for a considerably larger value of viscosity, namely $\mu=1$.  Here the solution does not steepen at all, but flattens out immediately after $t=0$, quickly resembling a slightly distorted Gaussian curve.  In the complex plane, where singularities for $\mu > 0$ must be simple poles, the effect of this regularisation is that for $\mu=0.1$ the closest pole to the real axis will initially move towards the axis but, unlike in the $\mu = 0$ case, subsequently will be redirected away from the real axis. On the other hand, for the larger value $\mu=1$, the closest pole to the real axis will immediately and forever move further away from the axis.  We will discuss this phenomenon in detail in this work.

The outline of the paper is as follows.  In section~\ref{sec:exact}, we illustrate the exact solution of (\ref{eq:maineq}), \eqref{eq:ic} in the complex plane with phase portraits and analytical landscapes. These results illustrate the singularity structure, which involves an array of poles in each quadrant that appears to continue indefinitely in the far field.  Section~\ref{sec:smalltime} is devoted to a rather comprehensive small-time analysis, which involves an inner region near the pole at $z=\mathrm{i}$ (and an identical reflection near $z=-\mathrm{i}$).  We are able to solve the inner problem exactly and show how the array of poles tends to line up on the anti-Stokes lines in the far field (of the inner problem).  In section~\ref{sec:poletrajectory} we focus on tracking the closest pole to the axis and make connections with the steepness of the solution on the real line.  We note the closest pole is typically the one of primary interest in applications as its distance to the real line determines the analyticity width of the solution \cite{caflisch2015complex, sulem1983}.  Section~\ref{sec:largetime} summarises the large-time solution behaviour, with a focus on the pattern of poles for late times.  Some of the analysis at the end of Section~\ref{sec:largetime} applies also for $t=\OO(1)$.
We are therefore able to provide a reasonably thorough understanding of the singularity dynamics across all time scales.  Section~\ref{sec:aaa} briefly explores how well a rational approximation works for this problem, using the AAA algorithm as an example.
Many of our calculations and numerical results throughout the paper are supported by further details in the Appendix and Supplementary Material. Finally, we close the paper in section~\ref{sec:discussion} with a summary of the key points and a discussion about future work. 
We mention issues that arise when the initial condition (\ref{eq:ic}) is replaced by one that has branch-point singularities and include a preliminary outline of the corresponding small-time asymptotics in the Appendix.

 All code to reproduce the numerical results in this work is available on  \href{https://github.com/DanielVandH/BurgersEquation.jl}{Github}.

\section{Exact solution with $u_0(x) = 1/(1+x^2)$}\label{sec:exact}

Burgers' equation \eqref{eq:maineq} can be linearised using the Cole-Hopf transformation \cite{cole1951,hopf1950}, transforming it into the heat equation. The solution is derived by setting $u = -2\mu v_x/v$, where $v(x, t)$ solves
\begin{equation}\label{eq:colehopf}
\pdv{v}{t}=\mu\pdv[2]{v}{x}, \quad \text{with}\quad v(x, 0) = \exp\left\{-\frac{1}{2\mu}\int_0^x u_0(\xi) \dxi\right\},~x\in\mathbb R,\,t>0.
\end{equation}
Using the Fourier transform to solve \eqref{eq:colehopf} for $v$ and rewriting in terms of $u$, we find \eqref{eq:exactsolution}. With the initial condition  \eqref{eq:ic}, equation \eqref{eq:exactsolution} becomes
\begin{equation}\label{eq:exactsolutionrunge}
u(z, t) = \dfrac{\displaystyle\int_{-\infty}^\infty \dfrac{z-s}{t}\exp\left\{-\dfrac{1}{2\mu}\left[\arctan s + \dfrac{(z-s)^2}{2t}\right]\right\} \ds}
{\displaystyle\int_{-\infty}^\infty \exp\left\{-\dfrac{1}{2\mu}\left[\arctan s + \dfrac{(z-s)^2}{2t}\right]\right\} \ds},
\end{equation}
where we now replace $x$ by $z$ to emphasise that this solution remains valid in the complex plane.

The integrals in \eqref{eq:exactsolutionrunge} are computed using a combination of Gauss-Hermite and Gauss-Legendre quadrature \autocite{burden2011numerical,sauer2012numerical} in \textsc{Julia} \autocite{bezanson2017julia} through the \texttt{FastGaussQuadrature.jl} package \autocite{townsend2015fast}.  A key step is to first make the change of variable $s=2\sqrt{\mu t}\bar{s}+z$ and then attempt to shift the contour on to the real $\bar{s}$-axis.  For $|\mathrm{Im}(z)| > 1$, the contour must be deformed around a branch point at $\bar{s}=(\i-z)/2\sqrt{\mu t}$.  Further details are provided in the Supplementary Material.  We note that the expressions derived in the Supplementary Material involve some highly oscillatory integrals, which could be evaluated using for example the methods presented by Dea{\~n}o et al. \autocite{deano2018computing}, although we naively apply Gauss-Legendre quadrature to such integrals; there will be some issues in our evaluation around the lines $\mathrm{Im}(z) = \pm 1$ (which turn out to be the Stokes lines in the limit $t\rightarrow 0^+$), although none too severe as we find throughout this paper that our numerical evaluations match the asymptotics of the solution extremely well.

\begin{figure}[h!]
\centering
\includegraphics[width=\textwidth]{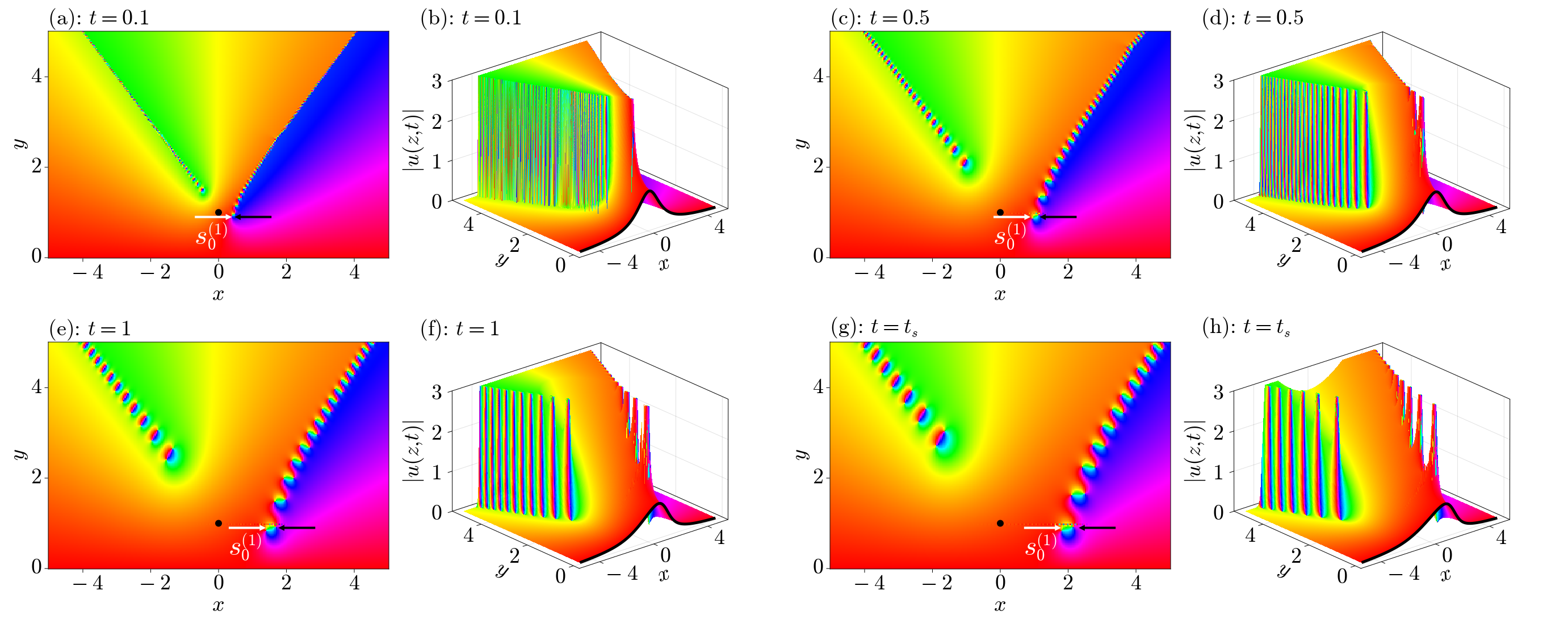}
\caption{Plots of the exact solution \eqref{eq:exactsolutionrunge} for $\mu = 0.1$, shown for times $t = 0.1$, $0.5$, $1$ and $t_s \approx 1.5396$.  For each time, the phase portrait is on the left. The closest pole to the real axis, $z=s^{(1)}_0$, and the associated zero are indicated by the white and black arrows, respectively. The point $z=\i$ is indicated by the solid black dot.  The analytical landscapes are on the right, together with black curves that indicate the solution on the real line. }\label{fig:complexsolutionplotviewsmallmu}
\end{figure}

\begin{figure}[h!]
\centering
\includegraphics[scale=0.2]{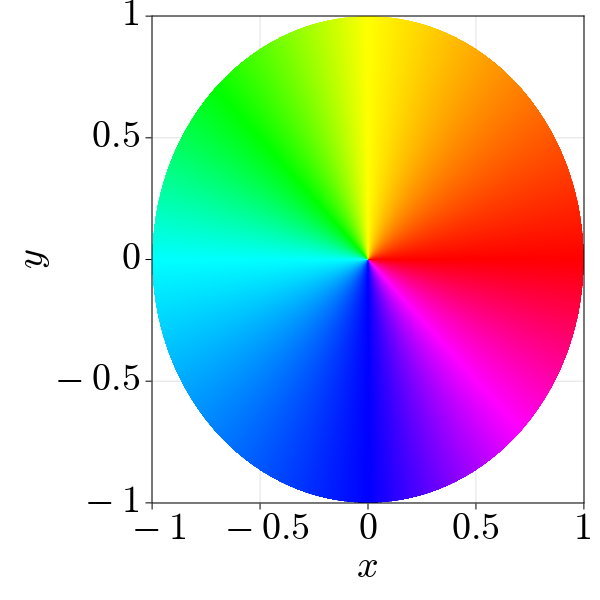}
\caption{Colour wheel used for plotting the phase portraits and analytical landscapes in the paper.}\label{fig:colourwheel}
\end{figure}

In figure~\ref{fig:complexsolutionplotviewsmallmu} we plot the solution \eqref{eq:exactsolutionrunge} for $\mu = 0.1$ at various times.  We restrict our attention here to the upper-half $z$-plane, remembering that the lower-half plane is a reflection about the real axis. For each time shown, in the left image we use phase portraits to visualise the solutions, with the phase of the solution at each point, $\arg(u(z, t))$, being used to colour each point \textcite{wegert2012visual,wegert2010phase}.  The colour wheel we use is given in figure~\ref{fig:colourwheel}. The right image for each time displays the analytical landscape of the solution, which uses the height at each point to represent the magnitude of the solution at the point.  For what follows, it is worth remembering that in a phase portrait a simple zero will appear locally like the colour wheel in figure~\ref{fig:colourwheel}, up to a rotation, while a simple pole will look similar except that the colours will be in reversed order as the pole is circumnavigated.

For each of the four times considered in figure~\ref{fig:complexsolutionplotviewsmallmu}, it is clear from closely inspecting the phase portraits and analytical landscapes that there is an array of poles in each of the first and second quadrants, and that every pole has a simple zero in close proximity (the zeros in the first quadrant lie immediately to the right of the poles, while the zeros in the second quadrant lie slightly below and to the left of the poles).  Our hypothesis is  that there are infinitely many poles (and zeros) in each array and that these all emerge from $z=\mathrm{i}$ at $t=0$ (an idea we explore in more detail in section~\ref{sec:smalltime}).  It seems clear from these images that all of these poles (and zeros) appear to be moving further away from $z=\mathrm{i}$ and, furthermore, become increasingly separated, as time increases.  Crucially, for this small value of viscosity, the closest singularity to the real axis, which we label $z=s^{(1)}_0(t)$ (indicated by the white arrow), sits in the first quadrant and appears to be moving closer to the real axis for these times.  The fourth time considered in this figure is $t= 8\sqrt 3/9$, which is the time at which blow-up would occur for $\mu=0$.  Here, for $\mu=0.1$, we see that this pole has not reached the real axis, as expected.  For larger times (not shown here), we find that this closest pole $z=s^{(1)}_0(t)$ changes its trajectory, beginning to move further away from the real axis.  We discuss this behaviour further in section~\ref{sec:poletrajectory}.

\begin{figure}[!h]
\centering
\includegraphics[width=\textwidth]{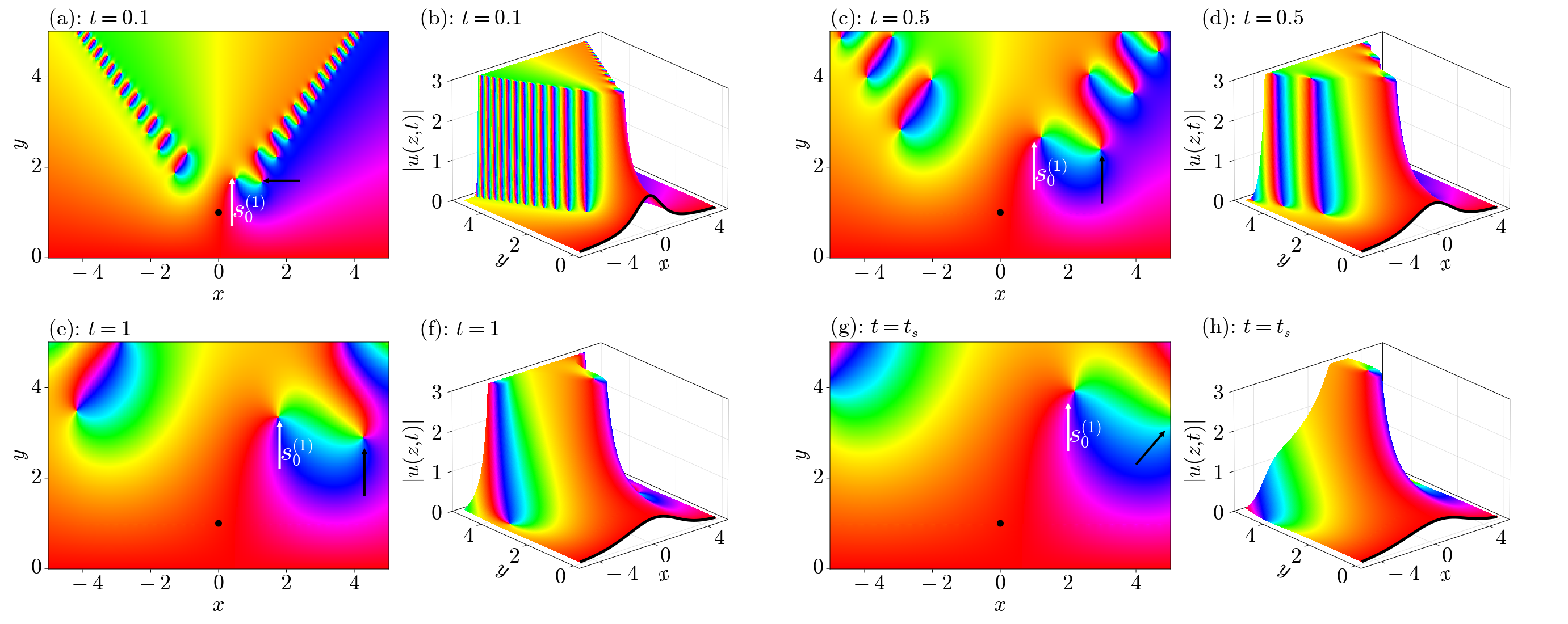}
\caption{Plots of the exact solution \eqref{eq:exactsolutionrunge} for $\mu = 1$, shown for times $t = 0.1$, $0.5$, $1$ and $t_s \approx 1.5396$.  For each time, the phase portrait is on the left.   The closest pole to the real axis, $z=s^{(1)}_0$, and the associated zero are indicated by the white and black arrows, respectively. The point $z=\i$ is indicated by the solid black dot.  The analytical landscapes are on the right, together with black curves that indicate the solution on the real line.}\label{fig:complexsolutionplotviewlargermu}
\end{figure}

The images in figure~\ref{fig:complexsolutionplotviewlargermu} show the solution \eqref{eq:exactsolutionrunge} for a larger value of viscosity, namely $\mu=1$.  Again, it is clear that, at the earliest time  $t=0.1$, there is an infinite array of simple poles (and associated simple zeros) in each of the first and second quadrants.  For the three subsequent times, the scale of the plots makes it more difficult to see the pattern, but again each of the poles appears to be moving away from $z=\mathrm{i}$ as time increases, although the separation distance between them is larger for $\mu=1$ when compared to $\mu=0.1$.  One significant difference from figure~\ref{fig:complexsolutionplotviewsmallmu} is that for $\mu=1$ the closest singularity to the real axis, $z=s^{(1)}_0(t)$ (indicated by the white arrow), appears to propagate so that it immediately moves further away from the real axis (i.e., it does not appear ever to move closer to the real axis).

In summary, these results suggest that at $t=0$ infinitely many simple poles are born at $z=\mathrm{i}$ and, for $t\ll 1$, appear to align themselves in a sequence with increasing distance from $z=\mathrm{i}$.  For a sufficiently small value of $\mu$, the closest pole to the real axis, $z=s^{(1)}_0(t)$, moves towards the axis and then eventually changes its direction and moves away from the real axis.  For a sufficiently large value of $\mu$, the closest pole to the real axis simply moves away from the real axis.  These observations suggest there is a borderline case (for $\mu$ somewhere between $0.1$ and $1$) in which the closest pole initially moves in the positive $z$-direction.  We address this issue in the following section.  Further, we explore how fast the poles are moving away from $z=\mathrm{i}$ in the small-time limit and what their precise pattern is, including their location and separation distance.

\section{Small time analysis, $t\ll 1$}\label{sec:smalltime}

In this section we provide a comprehensive analysis of the solution of \eqref{eq:maineq} with \eqref{eq:ic} in the small-time limit, focussing on the initial motion of the solutions' singularities.

\subsection{Naive expansion for the outer problem}
\label{sec:naiveexpansion}
To begin, we assume the naive expansion
\begin{equation}
\label{eq:naiveexpansion}
u(x, t) \sim u_0(x) + t u_1(x) + t^2u_2(x) + \ldots, \quad t \to 0^+.
\end{equation}
By substituting \eqref{eq:naiveexpansion} into the pde \eqref{eq:maineq}, we obtain
\begin{align*}
\left(u_1 + u_0u_0'\right) + \left(u_0'u_1 + u_0u_1' + 2u_2\right)t + \ldots &= \mu\left(u_0'' + tu_1'' + \ldots\right).
\end{align*}
Using our initial condition \eqref{eq:ic}, we find
\begin{equation}
\label{eq:naiveexpansionterms}
\begin{aligned}
u_1(x) &= \frac{2\left(-\mu + x + 3\mu x^2\right)}{\left(1+x^2\right)^3}, \\  u_2(x)& = \frac{60\mu^2 x^4 - 32\mu x + 48 \mu x^3 - x^2\left(120\mu^2 - 7\right) + 12\mu^2 - 1}{\left(1+x^2\right)^5}.
\end{aligned}
\end{equation}
The terms in $u_0$, $u_1$ and $u_2$ are the first few in what would be a divergent asymptotic expansion.  The divergence of this expansion is caused by the presence of singularities of the leading-order term (the initial condition $u_0(x)$) in the complex plane.  In our case, the singularities are at $z=\pm\mathrm{i}$.  Therefore, we expect an (optimally truncated) expansion (\ref{eq:naiveexpansion}) to apply not only on the real axis, but also as we move out to the complex plane, at least sufficiently far away from $z=\pm\mathrm{i}$ (we refine this claim below in subsection~\ref{sec:smalltimeWKB}).  Indeed, note that, as $z \to \i$,
\begin{align}
u_0 &\sim -\frac{\i}{2}\frac{1}{z-\i}+\frac14+\frac{\i}{8}(z-\i),
\label{eq:u0inner}\\
u_1 &\sim \left(-\frac14-\mu\i\right)\frac{1}{(z-\i)^3}-\frac{\i}{8}\frac{1}{(z-\i)^2}, \\
u_2 &\sim \left(\frac14-\frac{5\mu}{2}- 6\i\mu^2\right)\frac{1}{(z-\i)^5}.
\end{align}
As expected, the repeated double differentiation needed to calculate the terms in the expansion has led to singularities that are two orders higher in each case.
Therefore, the expansion (\ref{eq:naiveexpansion}) is no longer well ordered where $u_0 = \OO(tu_1)$, namely where
\begin{equation}
\label{eq:xidef}
\xi = \frac{z-\i}{t^{1/2}}=\OO(1).
\end{equation}
Hence, in the distinguished limit $z \to \i$, $t \to 0^+$,  we have
\begin{equation}
\label{eq:matchingconditionfromnaive}
u\sim \frac{1}{t^{1/2}}
\left[
-\frac{\i}{2\xi}+\left(-\frac14-\mu\i\right)\frac{1}{\xi^3}
+\left(\frac14-\frac{5\mu}{2}- 6\i\mu^2\right)\frac{1}{\xi^5}
+\ldots
\right]
+ \left[\frac14 - \frac{\i}{8\xi^2} + \cdots\right] + t^{1/2}\left[\frac{\i}{8}\xi + \cdots\right].
\end{equation}
Therefore, for $\xi = \OO(1)$, we write
\begin{equation}\label{eq:innerform}
u = \frac{1}{t^{1/2}}\Phi(\xi, t),
\end{equation}
and thus (\ref{eq:matchingconditionfromnaive}) provides far-field conditions for $\Phi$ as we discuss below.

Note that the exact solution \eqref{eq:exactsolutionrunge} can also be used to analyse \eqref{eq:maineq} in the limit $t \to 0^+$ using the method of steepest descents.  We very briefly summarise this approach in \ref{sec:smalltimesaddle} (and include further details in the Supplementary Material).

\subsection{Inner region $(z-\i)/t^{1/2}=\OO(1)$}
\label{sec:innerregion}

To continue our analysis, we consider the inner region $\xi=\OO(1)$, where $\xi$ is defined in (\ref{eq:xidef}).  We write \eqref{eq:innerform} so that \eqref{eq:maineq} can be rewritten exactly as
\begin{align}
t\pdv{\Phi}{t}-\frac12\Phi-\frac12\xi\pdv{\Phi}{\xi}+\Phi\pdv{\Phi}{\xi}=\mu\pdv[2]{\Phi}{\xi}.
\label{eq:innerpde}
\end{align}
For $t\ll 1$, we assume that
\[
\Phi(\xi, t) \sim \Phi_0(\xi) + t^{1/2}\Phi_1(\xi) + \ldots\quad\text{as}~t \to 0^+,
\]
where this $t^{1/2}$ scaling comes from (\ref{eq:matchingconditionfromnaive}).  From \eqref{eq:innerpde} we find
\begin{align}
\OO(1):~& -\sfrac{1}{2}\Phi_0 - \sfrac{1}{2}\xi\Phi_0' + \Phi_0\Phi_0' = \mu\Phi_0'', \label{eq:innerpdeconstantterm} \\
\OO(t^{1/2}):~& -\sfrac{1}{2}\xi\Phi_1' + \Phi_0\Phi_1' + \Phi_0'\Phi_1 = \mu \Phi_1'',\label{eq:innerpdelinearterm}
\end{align}
where primes indicate derivatives in $\xi$.  The far-field conditions, namely
\begin{equation}
\Phi_0\sim -\frac{\i}{2\xi}+\left(-\frac14-\mu\i\right)\frac{1}{\xi^3}
+\left(\frac14-\frac{5\mu}{2}- 6\mu^2\i\right)\frac{1}{\xi^5},
\quad
\Phi_1\sim \frac14 - \frac{\i}{8\xi^2},
\quad
\mbox{as}
\quad
\xi \to -\i\infty,
\label{eq:farfield}
\end{equation}
come from the inner limit (\ref{eq:matchingconditionfromnaive}) of the outer problem.  The limit is taken as $\xi \to -\i\infty$ so that the solution matches back onto the real line.  Note that, while we write out the problem here for $\Phi_1$, we only consider the leading-order solution for $\Phi_0$ in what follows.

\subsubsection{Exact solution to leading-order problem}

Integrating \eqref{eq:innerpdeconstantterm} once and enforcing the leading-order condition from (\ref{eq:farfield}), namely
\begin{equation}
\Phi_0\sim -\frac{\i}{2\xi}\quad\mbox{as}\quad \xi \to -\i\infty,
\label{eq:farfield2}
\end{equation}
gives
\begin{equation}\label{eq:firstorderpde}
\Phi_0^2 - \xi \Phi_0 = 2\mu \Phi_0' + \frac{\i}{2}.
\end{equation}
We note that \eqref{eq:firstorderpde} was also stated (but not analysed) in the appendix of Chapman et al.~\textcite{Chapman_2007}.  This Riccati equation can be solved exactly, again subject to (\ref{eq:farfield2}), to give \begin{equation}\label{eq:phisolutionparaboliccylinder}
\Phi_0(\xi) = \frac{1}{2\sqrt{2\mu}}\frac{U\left(\frac12-\frac{\i}{4\mu},\frac{\i\xi}{(2\mu)^{1/2}}\right)}{U\left(-\frac12-\frac{\i}{4\mu}, \frac{\i\xi}{(2\mu)^{1/2}}\right)},
\end{equation}
where $U$ is a parabolic cylinder function. Some details of this derivation are provided in \ref{eq:appendixinner}, along with a description of how we compute $U(a, z)$ for a given pair $(a, z)$ using the \texttt{HypergeometricFunctions.jl} package in \textsc{Julia} \autocite{Slevinsky2021}.  Further, we apply a Liouville-Green (WKB) argument to confirm in \ref{eq:appendixinner} that imposing (\ref{eq:farfield2}) represents two boundary conditions, as required.

\begin{figure}[h!]
\centering
\includegraphics[width=\textwidth]{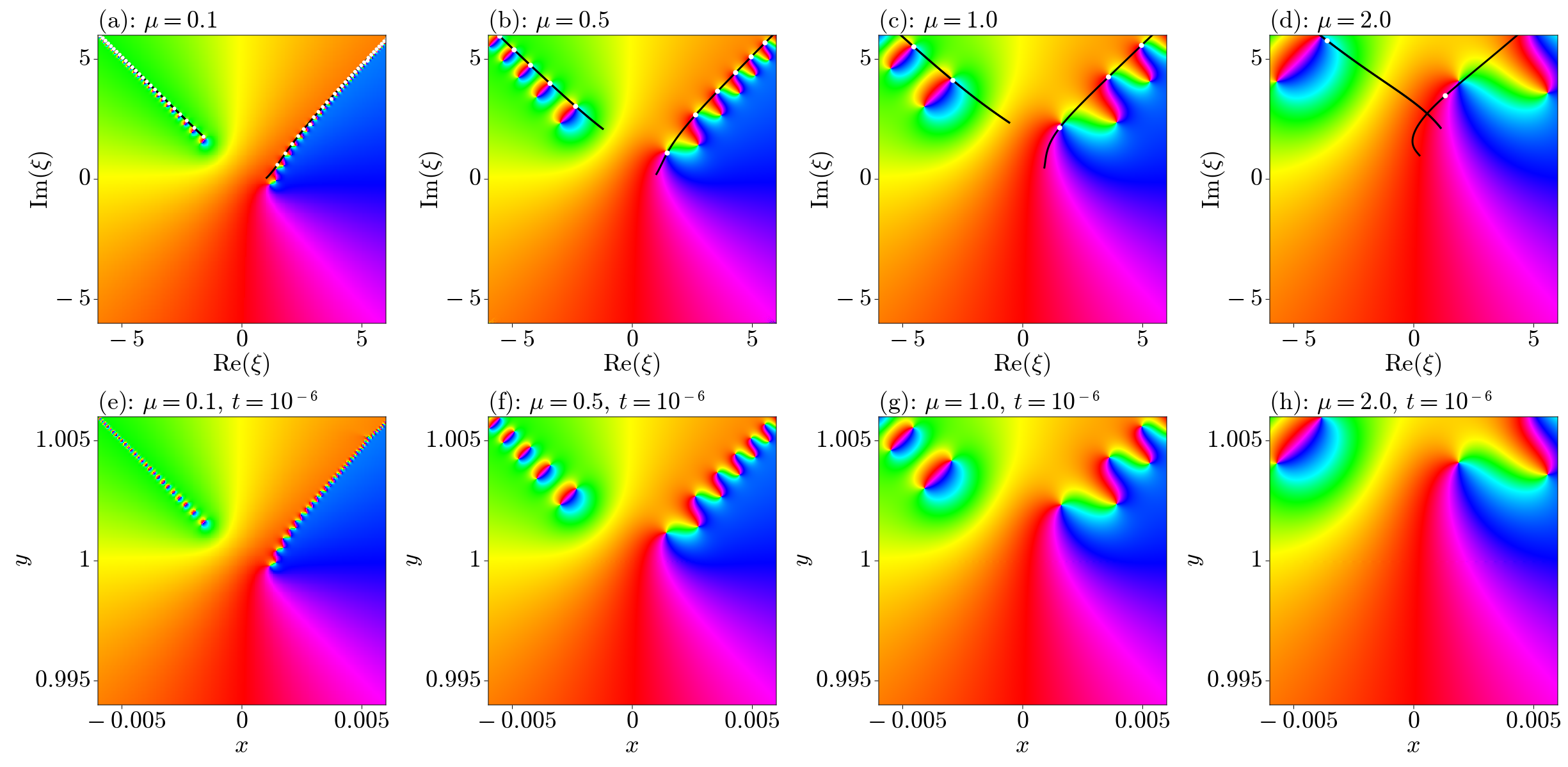}
\caption{(a)--(d) Phase portraits of $\Phi_0(\xi)$ in \eqref{eq:phisolutionparaboliccylinder} for $\mu=2$, $1$, $0.5$ and $0.1$. (i)--(l) Phase portraits of $u(z, t)$, zoomed in near $z=\i$, at the corresponding values of $\mu$ in (a)--(d) at $t = 10^{-6}$ for comparison. The black curves in the top row of plots show the path determined by the asymptotic angles \eqref{eq:vartheta} and \eqref{eq:vartheta2} for the positions of the poles for large $|\xi|$, and the white dots show the position of the poles determined by the transcendental equations \eqref{eq:transxi} and \eqref{eq:transxi2}.}\label{fig:phi0solutionparabolicgraphs}
\end{figure}

In figure~\ref{fig:phi0solutionparabolicgraphs}(a)--(d) we plot phase portraits of the inner solution \eqref{eq:phisolutionparaboliccylinder} for the four values
$\mu = 0.1$, $0.5$, $1$ and $2$.  Then in figure~\ref{fig:phi0solutionparabolicgraphs}(e)--(h) we show equivalent portraits for $u(z, t)$, via (\ref{eq:exactsolutionrunge}), which are drawn by zooming in close to $z=\i$ for a very small value of time, $t = 10^{-6}$.  We see that the solution for the inner problem is virtually indistinguishable from the exact solution for these values of $\mu$.  This comparison provides very strong evidence that the inner solution is correct. (We note that figure~\ref{fig:phi0solutionparabolicgraphs}(d) was computed with arbitrary precision, using the \texttt{ArbNumerics.jl} package in \textsc{Julia} which calls into the Arb C library \cite{sarnoff2018arb,johansson2017arb}, as the portrait could not be fully resolved using double precision.)

We make some preliminary observations of the inner solution plotted in figure~\ref{fig:phi0solutionparabolicgraphs}.  First, in both the first and second quadrants of the $\xi$-plane, there is an infinite array of poles that appear to approach rays at angles $\arg(\xi) = \pi/4$ and $\arg(\xi) = 3\pi/4$ as $|\xi|\rightarrow\infty$.
Second, we observe from these plots that, in (b)--(d), all poles appear in the upper-half plane, $\Im(\xi) > 0$.  Only in (a), for the case $\mu = 0.1$, do we see a pole in the lower-half $\xi$-plane.  Indeed, there appears to be a borderline value $\mu = \mu^* \approx 0.1468$, below which there is at least one pole in the lower-half $\xi$-plane.  We return to this point below.  Third, the colour spread around each pole can be observed to decrease with the value of viscosity $\mu$, supporting the fact that the strength of each pole depends linearly on $\mu$.

\subsubsection{Singularities of $\Phi_0$}

Concentrating on the poles we see in figure~\ref{fig:phi0solutionparabolicgraphs}, if we let $\xi_s$ be some pole of $\Phi_0$ and write $\Phi_0(\xi)\sim \phi_{-1}/(\xi-\xi_s) + \phi_0+\phi_1(\xi-\xi_s)$ as $\xi \to \xi_s$, then by substituting in (\ref{eq:innerpdeconstantterm}) we find
\begin{equation}\label{eq:finallocalpoleexpansion}
\Phi_0(\xi) \sim -\frac{2\mu}{\xi-\xi_s} + \frac12\xi_s + \left(\frac{8\mu - \xi_s^2 - 2\i}{16\mu}\right)\left(\xi-\xi_s\right),\quad\text{as}~\xi\to\xi_s.
\end{equation}
The analyses leading to figure~\ref{fig:phi0solutionparabolicgraphs} and the local expansion around simple poles $\xi_s$ in \eqref{eq:finallocalpoleexpansion} allow us to make the following conclusions.  Firstly, the strength of these simple poles (the modulus of the residue) is $2\mu$, implying that each pole's strength indeed grows in proportion with $\mu$ as was suggested from the colour spread in figure~\ref{fig:phi0solutionparabolicgraphs} around each pole. Secondly, the expansion \eqref{eq:finallocalpoleexpansion} verifies the existence of simple poles in the solution $u$ with the local behaviour (\ref{eq:localpole}); in the $z$-plane, we see that each pole propagates out from $z=\i$ on a trajectory asymptotic to $z = \i + t^{1/2}\xi_s$ as $t\rightarrow 0^+$.  This last point makes it clear why it is important to note that there is some borderline value $\mu = \mu^* \approx 0.1468$, below which there is at least one pole in the lower-half $\xi$-plane.  If we label the `lowest' of these, i.e. the one with the most negative imaginary part, $\xi=\xi_0$, then this implies that for $\mu<\mu^*$, the closest pole in the $z$-plane to the real axis, $z=s^{(1)}_0(t)$, initially moves towards the real axis as $s^{(1)}_0(t)\sim \i + t^{1/2}\xi_0$ as $t\rightarrow 0^+$. Indeed, as $\mu\rightarrow 0^+$, it appears that $\xi_0\rightarrow 1-\i$, so that $s^{(1)}_0(t)\sim \i + (1-\i)t^{1/2}$ in the double limit $\mu$, $t\rightarrow 0^+$, which matches the motion of the closest branch point for the inviscid problem \eqref{eq:inviscidsmall}.

This behaviour bears some resemblance  to that observed by Senouf \cite{senouf1997dynamicsA, senouf1997dynamics} for the initial condition $u(x, 0) = 4x^3-x/t_s$, whereby poles come in from infinity down the imaginary axis.  In that case, there exists some $\mu = \mu^*$ such that poles come in towards the real axis and then depart back up the imaginary axis after the shock time for $\mu < \mu^*$, but instead depart before the shock time when $\mu > \mu^*$.

A last remark is that we need to be careful about these interpretations for large $|\xi|$. We see that \eqref{eq:finallocalpoleexpansion}  holds for $|\xi_s| \gg 1$ in an increasingly small neighbourhood around $\xi = \xi_s$ as we move further into the far field. In this limit the term $1/(\xi-\xi_s)$ balances with $\xi_s$, thus \eqref{eq:finallocalpoleexpansion} becomes unreliable for $\xi = \xi_s + \OO(|\xi_s|^{-1})$.

\subsubsection{Large $\xi$ limit of $\Phi_0$}\label{sec:largxiPhi}

Recall that in order to match back onto the initial condition (\ref{eq:ic}), in the direction of the real $z$-axis, we enforce (\ref{eq:farfield2}).  The Liouville-Green (WKB) analysis results in \ref{eq:appendixinner} suggest that in the far field, an exponentially small term $C\,\xi^{-\i/2\mu}\,\e^{-\xi^2/4\mu}$ is switched on across the Stokes lines (the switching on of exponentially small terms across Stokes lines is a consequence of Stokes phenomenon and is described by Berry~\cite{berry1988} and many others, e.g.~\cite{oldedaalhuis1995}).  These Stokes lines occur where $\xi^2/4\mu$ is real and positive~\cite{dingle}, namely the positive and negative real $\xi$-axes.  Therefore the far-field behaviour $\Phi_0\sim -\i/2\xi$ as $\xi\rightarrow\infty$ will hold in all of the lower $\xi$-plane and then even further up to the rays $\arg(\xi)=\pi/4$ and $3\pi/4$, which represent anti-Stokes lines (where $\xi^2/4\mu$ is purely imaginary).

In our case we have an exact solution so therefore we can appeal to the far-field asymptotics of parabolic cylinder functions (in Section 12.9 of the Digital Library of Mathematical Functions~\cite{NIST:DLMF}, for example) to give
\begin{equation}
    \mbox{as}\,|\xi|\rightarrow\infty,\quad
    \Phi_0\sim\left\{
    \begin{array}{ll}
    \displaystyle -\frac{\i}{2\xi} & -5\pi/4<\mathrm{arg}(\xi)<\pi/4 \\
    \displaystyle \frac{2\sqrt{\mu\pi}\,\i\,\mathrm{e}^{\pi/8}}{\Gamma(-\i/4\mu)}\,\mathrm{e}^{-\i(|\xi|^2+2\ln|\xi|-\ln 2\mu)/4\mu}& \mathrm{arg}(\xi)=\pi/4\\
    \displaystyle \xi & \pi/4<\mathrm{arg}(\xi)<3\pi/4 \\
    \displaystyle
    \frac{2\sqrt{\mu\pi}\,\i\,\mathrm{e}^{-\pi/8}}{\Gamma(-\i/4\mu)}\,
    \mathrm{e}^{\i(|\xi|^2-2\ln|\xi|+\ln 2\mu)/4\mu} &
    \mathrm{arg}(\xi)=3\pi/4
    \end{array}
    \right.,
\label{eq:Phifarfield}
\end{equation}
where $\Gamma(z)$ is the gamma function.  The fact that $\Phi_0$ has different leading order asymptotic behaviours in different sectors of the $\xi$-plane is simply the Stokes phenomenon, and therefore not surprising.  However, the observation that $\Phi\sim \xi$ in part of the $\xi$-plane implies that the inner region does not match out to the initial condition in all directions.  From an outer perspective, this corresponds to $u\sim (z-\i)/t$ as $z\rightarrow\i$ for $t\ll 1$, $\pi/4<\arg(z-\i)<3\pi/4$ (cf. (\ref{eq:u0inner}), which holds for $z\rightarrow\i$, $t\ll 1$, $-5\pi/4<\arg(z-\i)<\pi/4$).  Further, a key conclusion is that the solution of (\ref{eq:maineq}) with (\ref{eq:ic}) is not asymptotic to $1/(1+z^2)$ everywhere in the $z$-plane as $t\rightarrow 0^+$.

In order to approximate the location of the poles of $\Phi_0$, we look for zeros of the denominator of (\ref{eq:phisolutionparaboliccylinder}).  Far-field properties of parabolic cylinder functions give \cite{NIST:DLMF}
\begin{align*}
U\left(-\frac12-\frac{\i}{4\mu}, \frac{\i\xi}{(2\mu)^{1/2}}\right)
&\sim \left(\frac{\i\xi}{\sqrt{2\mu}}\right)^{\i/4\mu}\e^{\xi^2/8\mu}
\left[1+\frac{\sqrt{2\pi}\,\i}{\Gamma(-\i/4\mu)}
\left(\frac{\xi}{\sqrt{2\mu}}\right)^{-1-\i/2\mu}\e^{-\xi^2/4\mu}
\right],
\quad\xi\rightarrow\infty\quad \mbox{(first quadrant).}
\end{align*}
If we let the term in the square brackets be $G$, then to locate poles of $\Phi_0$ in the far field near $\arg(\xi)=\pi/4$, we need to solve $G=0$ asymptotically.

The term $\Gamma(-\i/4\mu)$ is complex whose real and imaginary part must be determined numerically.  However, using the result $|\Gamma(\i y)|^2=\pi/(y\sinh \pi y)$, we can identify the modulus exactly, so that
$$
\frac{1}{\Gamma(-\i/4\mu)}=\left(
\frac{\sinh(\pi/4\mu)}{4\mu\pi}
\right)^{1/2}\e^{-\i\alpha},
$$
where $\alpha$ is real (with $\alpha\rightarrow \pi^+/2$ as $\mu\rightarrow\infty$).
We write $\xi=|\xi|\e^{\i\vartheta}$, then, starting with the first quadrant,
 \begin{align*}
G&=1+\frac{\sinh^{1/2}(\pi/4\mu)}{|\xi|}
\e^{(-|\xi|^2\cos 2\vartheta+2\vartheta)/4\mu}
\left[\cos\left(
\frac{|\xi|^2\sin 2\vartheta+2\ln|\xi|-\ln(2\mu)}{4\mu}+\vartheta-\frac{\pi}{2}+\alpha
\right)\right.
\\
 &- \left.\i \sin\left(
\frac{|\xi|^2\sin 2\vartheta+2\ln|\xi|-\ln(2\mu)}{4\mu}+\vartheta-\frac{\pi}{2}+\alpha
\right)\right]
\end{align*}
For the imaginary part of $G$ to vanish, we have
\begin{equation}
\sin\left(
\frac{|\xi|^2\sin 2\vartheta+2\ln|\xi|-\ln(2\mu)}{4\mu}+\vartheta-\frac{\pi}{2}+\alpha
\right)=o(1),
\label{eq:cossin1}
\end{equation}
\begin{equation}
\cos\left(
\frac{|\xi|^2\sin 2\vartheta+2\ln|\xi|-\ln(2\mu)}{4\mu}+\vartheta-\frac{\pi}{2}+\alpha
\right)\sim -1,
\quad|\xi|\rightarrow\infty,
\label{eq:cossin2}
\end{equation}
and, therefore, for the real part of $G$ to vanish we need
$$
1-\frac{\sinh^{1/2}(\pi/4\mu)}{|\xi|}
\e^{(-|\xi|^2\cos 2\vartheta+2\vartheta)/4\mu}=o(1),
$$
or
\begin{equation}
\vartheta\sim \frac{\pi}{4}+\frac{2\mu}{|\xi|^2+1}
\left(\ln|\xi|-\frac{\pi}{8\mu}-\frac12\ln(\sinh(\pi/4\mu))
\right)
\quad\mbox{as}\quad|\xi|\rightarrow\infty,
\label{eq:vartheta}
\end{equation}
where terms neglected are $o(|\xi|^{-2})$.
Poles of $\Phi_0$ near $\vartheta=\pi/4$ lie approximately on this path.
Combining with (\ref{eq:cossin1})-(\ref{eq:cossin2}) gives the transcendental equation
\begin{equation}
-\frac{|\xi|^2+1}{2\mu}\tan
\left(
\frac{|\xi|^2+2\ln|\xi|-\ln(2\mu)}{4\mu}-\frac{\pi}{4}+\alpha
\right)
=\ln|\xi|-\frac{\pi}{8\mu}-\frac12\ln(\sinh(\pi/4\mu)),
\label{eq:transxi}
\end{equation}
which provides a means to predict the location of the poles along (\ref{eq:vartheta}).  The very first approximation to this equation gives $|\xi|^2 \sim 8n\mu\pi$, where $n$ is an integer, which leads to the asymptotic spacing $\sqrt{2\mu\pi/n}$ as $n\rightarrow\infty$.

A similar analysis for the poles in the second quadrant, where the poles tend to approach the anti-Stokes line $\vartheta = 3\ppi/4$, leads to the expression for $\vartheta$ in the second quadrant,
\begin{align}\label{eq:vartheta2}
\vartheta \sim \frac{3\ppi}{4} + \frac{2\mu}{|\xi|^2 - 1}\left(-\ln|\xi| - \frac{\ppi}{8\mu} + \frac12\ln(\sinh(\ppi/4\mu))\right), \quad |\xi| \to \infty,
\end{align}
where we have again neglected terms $o(|\xi|^{-2})$. Similarly, the modulus $|\xi|$ of each pole in the second quadrant is approximately the solution to the following transcendental equation,
\begin{align}\label{eq:transxi2}
\frac{|\xi|^2-1}{2\mu}\tan\left(\frac{2\ln|\xi| - |\xi|^2 - \ln(2\mu)}{4\mu} + \frac{\ppi}{4} + \alpha\right) = \ln|\xi| - \frac{\ppi}{8\mu} - \frac12\ln(\sinh(\ppi/4\mu)).
\end{align}
The spacing between poles here is also $\sqrt{2\mu\ppi/n}$ as $n \to\infty$.

We show the results of this analysis in the first and second quadrants of the $\xi$-plane in Figure \ref{fig:phi0solutionparabolicgraphs}(a)-(d), using black curves to show the paths \eqref{eq:vartheta} and \eqref{eq:vartheta2}, and white dots to indicate solutions to the transcendental equations \eqref{eq:transxi} and \eqref{eq:transxi2}. Keeping in mind that these results are asymptotic in the limit $|\xi|\rightarrow\infty$, we see that the predictions are remarkably accurate, capturing the pole positions very well. This is particularly true for $\mu=0.5$ and $1$, where predictions appear extremely good on the scale of the figure.

\subsection{Liouville-Green (WKB) for small time}\label{sec:smalltimeWKB}

Very briefly, we outline a crude Liouville-Green (WKB) analysis for small time in order to highlight exponentially small contributions that are switched on across Stokes lines.  Away from $z=\pm i$, we can linearise around the initial condition by writing
$$
u\sim \frac{1}{1+z^2}+U(z,t),
$$
where our Liouville-Green ansatz is
$$
U\sim \e^{-W_0(z)/t+W_1(z)}.
$$
After substituting into Burgers' equation and gathering terms of order $t^{-2}$ and $t^{-1}$, we find
$$
W_0=\frac{(z\mp\i)^2}{4\mu},
\quad
W_1=-\frac{1}{2}\log(z\mp\i)+\frac{1}{2\mu}\arctan z +\,\mathrm{constant}.
$$
There will be two possible terms, each associated with either the plus or minus signs.

Sticking to the upper-half plane, we therefore see that, away from $z=\i$, an exponentially small term of the form
$$
K(z-\i)^{-1/2}\e^{(\arctan z)/2\mu}\e^{-(z-\i)^2/4\mu t}
$$
switches on across Stokes lines $\Im(z)=1$ in the limit $t\rightarrow 0^+$.  These Stokes lines come from setting the singulant $(z-\i)^2/4\mu$ to be real and positive.
We note that this analysis suggests that we cannot expect the naive outer expansion (\ref{eq:naiveexpansion}) to apply past the anti-Stokes lines (when the singulant $(z-\i)^2/4\mu$ is imaginary), namely when $\pi/4<\arg(z-\i)<3\pi/4$.  This observation is consistent with our conclusion in subsection~\ref{sec:largxiPhi} that our solution of (\ref{eq:maineq}) and (\ref{eq:ic}) does not approach the initial condition $1/(1+z^2)$ in all parts of the $z$-plane as $t\rightarrow 0^+$.

\subsection{Summary of small-time behaviour}

In summary, a naive series expansion (\ref{eq:naiveexpansion}) in powers of $t$ breaks down near singularities of the leading-order term, namely at $z=\pm\i$.  We focus here on the upper plane only, due to symmetry.  By expanding out each term near $z=\i$, we arrive at scalings for an inner region (\ref{eq:innerform}) and derive matching conditions in the far field of the inner variable $\xi$ defined in (\ref{eq:xidef}).  The leading-order inner problem for $\Phi_0$ (also posed in the appendix of Chapman et al.~\textcite{Chapman_2007}) has an exact solution (\ref{eq:phisolutionparaboliccylinder}) in terms of parabolic cylinder functions, which has an array of simple poles in each of the first and second quadrants of the $\xi$-plane that tend to lie on the anti-Stokes lines $\arg(\xi)=\pi/4$ and $\arg(\xi)=3\pi/4$ as $|\xi|\rightarrow\infty$.  Using known asymptotic results about parabolic cylinder functions, we are able to explore the behaviour of $\Phi_0$ in the far field, including an approximate location of the poles.

A key point to emphasise is that infinitely many poles of $u(z,t)$ emerge spontaneously from $z=\i$ at $t=0^+$, despite there being only one pole at $t=0$.  The poles all initially move from $z=\i$ with speed $|\xi_s|/2t^{1/2}$ and direction $\arg(\xi_s)$, where $\xi_s$ is the corresponding pole of the inner problem.  The spacing between the poles is roughly $\sqrt{2\mu\pi/n}\,t^{1/2}$, which becomes increasingly small as viscosity $\mu$ decreases.  Therefore, this small-time analysis is able to provide a rather comprehensive summary of how poles are created in our problem and their early behaviour.

A further important point is that, from the images in figures~\ref{fig:complexsolutionplotviewsmallmu} and \ref{fig:complexsolutionplotviewlargermu} and the far-field expansion (\ref{eq:Phifarfield}) of the inner problem for $\Phi_0$, the analytic continuation of the solution of Burgers' equation (\ref{eq:maineq}) with (\ref{eq:ic}) is an analytic function in the $z$-plane that does not approach $1/(1+z^2)$ everywhere in the $z$-plane as $t\rightarrow 0^+$.  Instead, $u$ grows like $u\sim z/t$ in the limits $t\rightarrow 0^+$ and $z\rightarrow\infty$ in the wedge $\pi/4<\arg(z)<3\pi/4$, which again does not match with $1/(1+z^2)$ as $t\rightarrow 0^+$.

\section{Trajectory of nearest singularity to the real axis}\label{sec:poletrajectory}

In this section, we focus on the trajectory of the pole closest to the real axis, $z=s^{(1)}_0(t)$.  In subsection~\ref{sec:methodsteepestdescents}, we summarise how we do this using the method of steepest descents for the regime $\mu\ll 1$.  In the subsequent two subsections we make connections between the trajectory of $z=s^{(1)}_0(t)$ in the complex plane and the steepness of the solution profile on the real line.

\subsection{Method of steepest descents for $\mu\ll 1$}
\label{sec:methodsteepestdescents}

Here we apply steepest descents (the saddle-point method) in the spirit of Senouf~\cite{senouf1997dynamicsA} and Weideman~\cite{weideman2021dynamics}, for example.  A difference is that the poles all lie on the imaginary axis for their initial conditions, simplifying the geometry somewhat.  Keeping in mind we wish here to track singularities of $u$, rather than the full solution (\ref{eq:exactsolutionrunge}), we shall be focussing on zeros of the denominator
\begin{equation}
\label{eq:denominatordmudependence}
D(z,t;\mu) = \int_{-\infty}^\infty \e^{h(s)/\mu} \ds, \quad\text{where}\quad h(s) = -\frac12\arctan(s)-\frac{(z-s)^2}{4t}.
\end{equation}
For $\mu \ll 1$, the saddle points of $h(s)$ come from solving $h'(s) = 0$, leading to the cubic equation
\begin{equation}
\label{eq:cubicequationsaddlepoints}
s^3  - zs^2 + s + t - z = 0.
\end{equation}
For $z \in \mathbb R$, $t < 8\sqrt 3/9$, one saddle point is real and the other two complex conjugates.  Then, for $z \in \mathbb R$, $t > 8\sqrt 3/9$, between the two caustics there are three distinct real saddle points.  Regardless,  there are no singularities in $u$ for $z \in \mathbb R$ and so we concentrate on the complex plane.

For each $(z,t)$, with $z \in \mathbb C$ and $t>0$, there are, in general, three complex solutions of (\ref{eq:cubicequationsaddlepoints}), which we label $s_j \in \mathbb C$, and so three saddle-point contributions.
Our method for obtaining the relevant contributions will involve considering the local behaviour of $h(s)$ near each saddle point to determine the angle at which the contour should cross each saddle \autocite[Section 6.6]{bender1999advanced}.  This incident angle $\delta_j$ allows us to parameterise a tangent line to the contour $s - s_j \sim \sigma \e^{\i\delta_j}$ as $s \to s_j$, where $\sigma = |s-s_j|$.  If we write $h''(s_j) = |h''(s_j)| \e^{\i\alpha_j}$ then, using $h'(s_j) = 0$,
\begin{equation}\label{eq:saddlethetacondition}
h(s) - h(s_j) \sim \frac{1}{2}|h''(s_j)|\sigma_j^2\,\e^{\i(2\delta_j+\alpha_j)}.
\end{equation}
Since the contour is a steepest-descent contour, the imaginary part of $h(s)$ is constant, equalling $\Im h(s_j)$, and the real part is decreasing away from $s_j$. Therefore, it must be that $\sin(2\delta_j+\alpha_j) = 0$ and $\cos(2\delta_j+\alpha_j) < 0$. This gives two solutions for $\delta_j$,
\begin{equation}\label{eq:saddlepointangles}
\delta_j^{(1)} = -\frac{\alpha_j}{2}+\frac{\pi}{2}\quad\text{and}\quad \delta_j^{(2)} = -\frac{\alpha_j}{2} + \frac{3\pi}{2}.
\end{equation}
These values of $\delta_j$ are always in opposite quadrants. If we write the tangent line as $s = s_j + \varepsilon\,\e^{\i\delta_j^{(k)}}$, then selecting the value $\delta_j^{(k)} \in \{\delta_j^{(1)}, \delta_j^{(2)}\}$ such that $\cos(\delta_j^{(k)}) > 0$ will allow us to take $\varepsilon \in \mathbb R$ such that $\varepsilon > 0$ corresponds to the appropriate side of the tangent line.

In summary, to locate a zero of $D$ in (\ref{eq:denominatordmudependence}) we must look for values of $z$ for which we have two saddle-point contributions in the $s$-plane of the same size, so they can cancel to leading order.  For there to be two equally dominant saddle-point contributions, then we require the real parts of $h$ evaluated at each saddle point to be the same.  Numerically, this cancellation does not occur exactly.  Instead, for a given $z$, we compute $s_i$ for $i=1,2,3$  directly from \eqref{eq:cubicequationsaddlepoints}, and then we find the pair of distinct indices $(j_1, j_2)$ that correspond to the largest two contributions (i.e., the pair that minimises $|\mathrm{Re}(h(s_{j_1}) )- \mathrm{Re}( h(s_{j_2}))|$).  These indices define the saddle points that are used.  Using the two saddle points $s_{j_1}$ and $s_{j_2}$, following the usual method of steepest descents, the asymptotic expansion of the denominator \eqref{eq:denominatordmudependence} at a pole $z$ of $u$  is given by
\begin{equation}\label{eq:saddlepointasymptoticexpansionsmallmu}
D(\mu) \sim \sqrt{2\mu\pi}\sum_{j=1}^2 |h''(s_{j_i})|^{-1/2}\e^{\i\delta_{j_i} + h(s_{j_i})/\mu},\quad \text{as}~\mu\to 0^+,
\end{equation}
where $\delta_{j_i}$ is the value of $\delta_j$ in \eqref{eq:saddlepointangles} such that $\cos(\delta_{j_i}) > 0$. Further details on computing \eqref{eq:saddlepointasymptoticexpansionsmallmu} are given in \ref{sec:appendixsteepest}. When it comes to implementing this strategy computationally, we have to worry only about the selection of $s_{j_i}$ and $\delta_{j_i}$ for $i=1,2$. We compute $\arg h''(s_{j_1})$ and $\arg h''(s_{j_2})$ so that we can select the angle in \eqref{eq:saddlepointangles} that has a positive cosine, thus allowing us to compute \eqref{eq:saddlepointasymptoticexpansionsmallmu}.

Now we can apply the above strategy to track the closest pole in the complex plane to the real line. We start by obtaining a reliable estimate of the closest pole at a large time, say $t=2$ (which we may obtain by visually inspecting the numerical evaluation of \eqref{eq:exactsolutionrunge}). Then at a time $t_j = t - j\Delta t$, $j =1,2,\ldots$, we use a two-point linesearch on $D(z, t_j) = 0$, where we now approximate $D$ using \eqref{eq:saddlepointasymptoticexpansionsmallmu}, where the initial guess for the pole location $z$ is based on the estimated pole at time $t_{j-1}$. This procedure is continued until the solution breaks down due to ill-conditioning for small time. The gradients required in this linesearch are computed using automatic differentiation with the \texttt{ForwardDiff.jl} package in \textsc{Julia} \autocite{revels2016forward}.

\subsection{Role of borderline value $\mu=\mu^*$}

Recall that our small-time analysis suggests that $s^{(1)}_0\sim \i+t^{1/2}\xi_0$, which provides two pieces of information.  First, the closest pole to the real axis initially moves in a straight line in the $z$-plane from $z=\i$ in the direction $\mathrm{arg}(\xi_0)$.  In the borderline case $\mu=\mu^*$, the pole $\xi_0$ lies on the positive real $\xi$-axis and so the direction of $z=s^{(1)}_0(t)$ is initially parallel to the $x$-axis, moving in the positive $x$-direction.  Second, we see the distance that $z=s^{(1)}_0(t)$ takes from the $x$-axis evolves like $1+\Im(\xi_0)\, t^{1/2}$ as $t\rightarrow 0^+$.  In the borderline case $\mu=\mu^*$, clearly this distance is $1+\OO(t)$ in the limit.

\begin{figure}[h!]
\centering
\includegraphics[width=\textwidth]{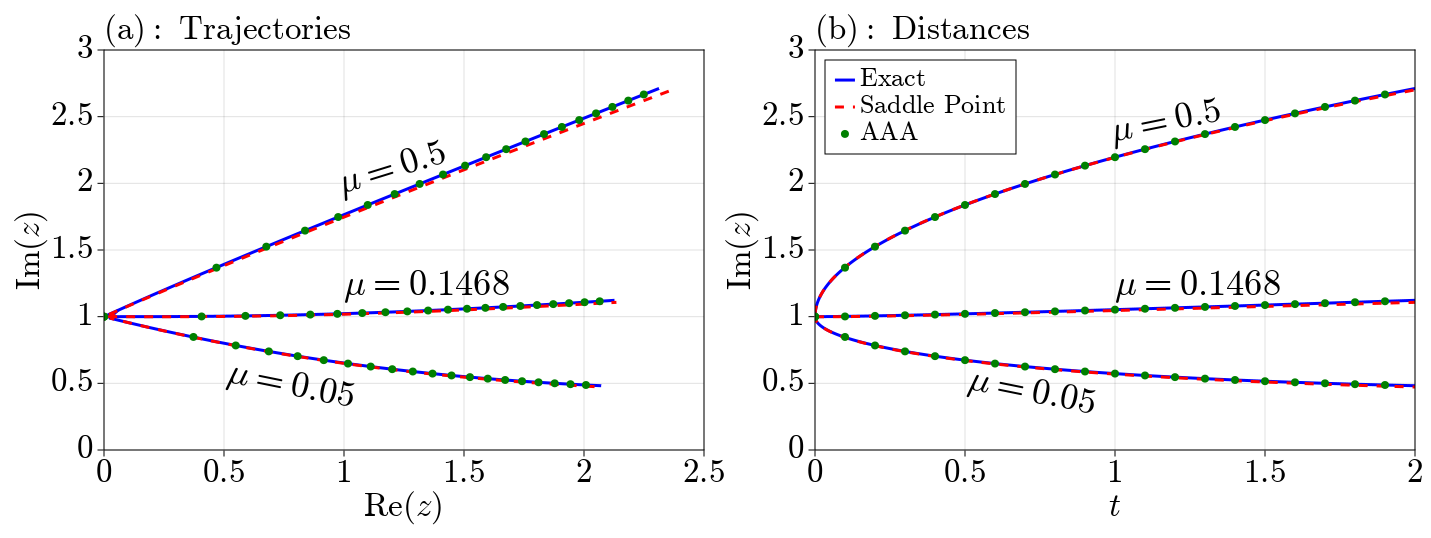}
\caption{Tracking the closest pole to the real line of $u(z,t)$, which we label $z=s^{(1)}_0$.  The ``exact'' trajectory (blue solid) is computed by tracking the zeros of the denominator in \eqref{eq:exactsolutionrunge}, computed using Gauss-Legendre quadrature.  The saddle point trajectory (red dashed) is computed using the saddle point approximation \eqref{eq:saddlepointasymptoticexpansionsmallmu}. The AAA results (green dots) are computed using the algorithm described in section \ref{sec:aaa}.  In the left plot the right-most point corresponds to $t=2$.}\label{fig:saddlepointtrajectorytracked}
\end{figure}

We observe this behaviour in figure~\ref{fig:saddlepointtrajectorytracked} for three solutions of (\ref{eq:maineq}) with \eqref{eq:ic} up to $t=2$, which is chosen to be some representative $\OO(1)$ time. These results, determined by computing roots of the denominator $D(z,t;\mu)$ of (\ref{eq:exactsolutionrunge}), are drawn as (blue) solid curves.  For the small value $\mu=0.05$, the closest pole $z=s^{(1)}_0(t)$ clearly evolves towards the real axis over this time scale, while for the larger value $\mu=0.5$ it initially moves away.  In the borderline case $\mu=\mu^*\approx 0.1468$, the trajectory appears initially horizontal, as expected, although it tends to slope slightly upwards as time increases.  Note that these trajectories are similar to those for the initial conditions $u(x, 0) = 4x^3 - x/t_s$ \cite{bessis1984pole,bessis1990complex,senouf1997dynamicsA,senouf1997dynamics} or $u(x, 0) = -\sin(x)$ \cite{sulem1983,weideman2021dynamics}, in the sense that the closest pole initially moves in along the imaginary axis up to some minimum distance before changing direction and going back up along the imaginary axis.  The difference is that for our problem the trajectory of the closest pole is two dimensional (i.e., not restricted to a linear path).

Also included in figure~\ref{fig:saddlepointtrajectorytracked} as (red) dashed curves are results calculated using the method of steepest descents.  When comparing with the (blue) solid curves which come from the exact solution (\ref{eq:exactsolutionrunge}), these steepest-descent predictions are indeed very good, even for $\mu=0.5$ which is not particularly small.  The data in the form of (green) dots comes from a rational approximation via the AAA algorithm.  We return to these data in section~\ref{sec:aaa}.

\begin{figure}[h!]
\centering
\includegraphics[width=\textwidth]{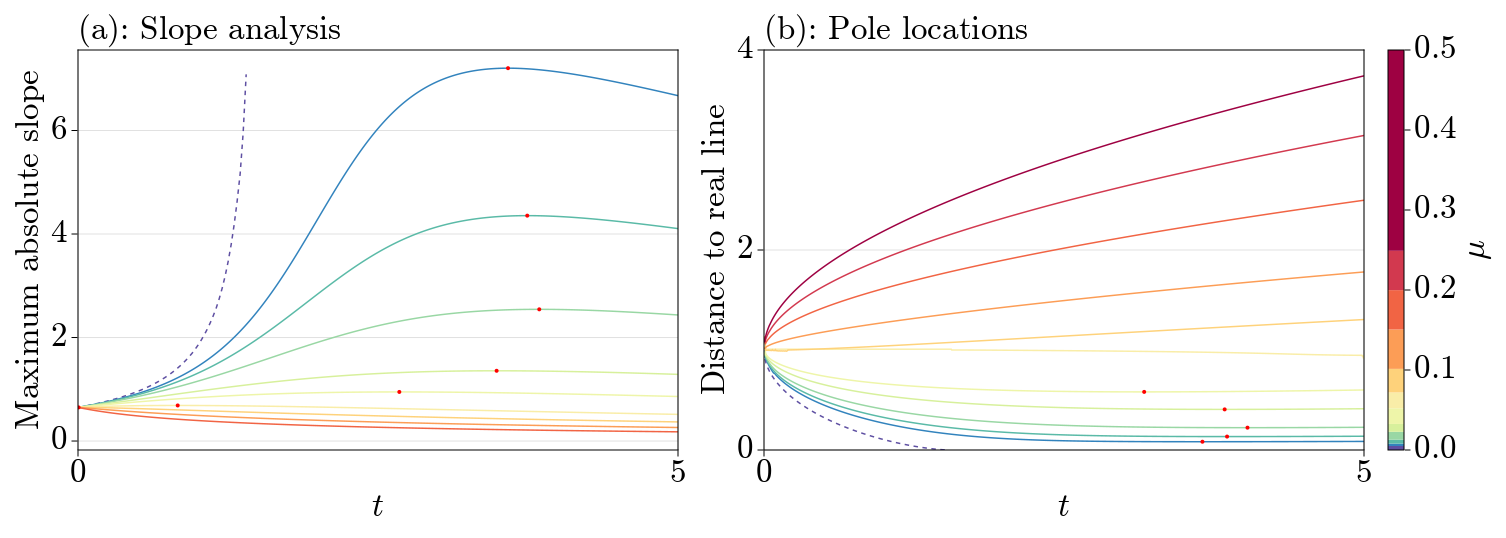}
\caption{(a) Maximum absolute slope for the solution $u(x,t)$ on the real line versus time, computed for a variety of values of $\mu$ as per the colour bar on the far right.  The (red) solid dots indicate local maxima on the individual curves for each $\mu$.  The (blue) dashed curve is for $\mu=0$, which continues to increase until $t_s \approx 1.5396$.  (b) The distance of the closest pole $z=s^{(1)}_0$ to the real axis versus time, again for values of $\mu$ as per the colour bar.  The (red) solid dots indicate local minima on the individual curves for each $\mu$. The (blue) dashed curve is for $\mu=0$, which intersects the $t$-axis at $t=t_s$.}
\label{fig:advectionslopepoleanalysis}
\end{figure}

We are interested in relating these properties to behaviour of the solution on the real line.  In particular, we explore the possible connections between the closest pole to the real axis, $z=s^{(1)}_0(t)$, and the steepness of the solution on the real line.  Recall, from figure~\ref{fig:compareinviscidviscous}, that the solution profile for $\mu=0.1$ initially begins to steepen (before later flattening out) while the profile for $\mu=1$ flattens out immediately.  This behaviour is summarised in figure~\ref{fig:advectionslopepoleanalysis} for a broader spread of viscosity values.  In (a), we plot the maximum (absolute value) of the slope on the real solution versus time.  These are determined by first computing the slope of the solution using automatic differentiation with the \texttt{ForwardDiff.jl} package in \textsc{Julia} \autocite{revels2016forward}, and then computing the maximum slopes using \texttt{Optim.jl} \autocite{mogensen2018optim}.  Clearly for $\mu=0$ (blue dashed curve) the maximum slope increases indefinitely until the inviscid blow-up time $t=t_s=8\sqrt{3}/9$, while for small values of $\mu$ the maximum slope increases in time (before eventually decreasing).  On the other hand, for larger values of $\mu$ the slope simply decreases in time.  Clearly there is a borderline case $\mu=\tilde{\mu}$ which divides these two different qualitative behaviours.  Studying this behaviour of the real solution numerically, we estimate the value to be $\tilde{\mu}\approx 0.1458$.  The numerical similarity between $\mu^*\approx 0.1468$ and $\tilde{\mu}\approx 0.1458$ is quite remarkable, suggesting a very strong link between the steepness of the real solution at small times and the initial trajectory of the closest singularity to the real line.  The existence of this link is perhaps not surprising, but certainly the closeness of these numerical values is interesting.

We remark briefly that an alternative measure of solution regularity on the real line is {\em enstrophy}, defined by
$$
E(t)=\frac{1}{2}\int_{-\infty}^\infty
\left(\frac{\partial u}{\partial x}\right)^2\,\mathrm{d}x
$$
\cite{ludoering2008} (noting that some definitions leave off the $1/2$ out the front), with the property
$$
\frac{\mathrm{d}E}{\mathrm{d}t}=-\mu \int_{-\infty}^\infty \left(\frac{\partial^2 u}{\partial x^2}\right)^2\,\mathrm{d}x
-\frac{1}{2}\int_{-\infty}^\infty \left(\frac{\partial u}{\partial x}\right)^3\,\mathrm{d}x.
$$
Solutions for $u(x,t)$ are smooth as long as the enstrophy is finite.  Researchers have used Burgers' equation to determine bounds on enstrophy and to explore the sharpness of these estimates \cite{ayala2011,pelinovsky2012}.  Links with the closest singularity in the complex plane are made in \cite{ayala2011}.
In this way, Burgers' equation plays the role of a much-simplified analogue of the three-dimensional Navier-Stokes equations, for which the regularity problem has attracted considerable interest, e.g. \cite{protas2022}.  In our case, we can produce plots of enstrophy versus time that appear qualitatively similar to figure~\ref{fig:advectionslopepoleanalysis}(a).  That is, when $\mu=0$, the enstrophy continues to increase without bound until the shock forms at $t=t_s$.
For sufficiently small $\mu>0$, the enstrophy increases and then decreases in time, while for sufficiently large $\mu$ the enstrophy simply decreases in time.
We have included the associated figure in the Supplementary Material.

\subsection{Balancing advection and diffusion}

We pursue these ideas further.  For $\mu<\tilde{\mu}$, the maximum steepness of the solution on the real line initially increases in time, reaches a maximum, and then decreases.  That is, as is well known for solutions to Burgers' equation, even if wave steepening via advection dominates for early time, ultimately smoothing via diffusion becomes significant.  We see the turning points of this behaviour in figure~\ref{fig:advectionslopepoleanalysis}(a) as red dots.  For example, for the five values $\mu=0.01$, $0.01590$, $0.02574$, $0.04542$ and $0.06510$, the turning points occur at times $t=3.5836$, $3.7437$, $3.8438$, $3.4885$ and $2.6777$, respectively.  Turning to figure~\ref{fig:advectionslopepoleanalysis}(b), the distance of the closest $z=s^{(1)}_0(t)$ to the real axis is plotted versus time for the same values of $\mu$ as in (a).  Here we see for $\mu>\mu^*$, each of these curves decreases in time, reaches a minimum, and then increases, with the minimum value indicated by the red dot.  For the five values $\mu=0.01$, $0.01590$, $0.02574$, $0.04542$ and $0.06510$, the turning points are estimated to be at $t=3.6537$, $3.8589$, $4.0290$, $3.8388$ and $3.1682$, respectively.  Thus we see that while there is a qualitative correspondence between the location of the local maxima in figure~\ref{fig:advectionslopepoleanalysis}(a) and the local minima in figure~\ref{fig:advectionslopepoleanalysis}(b), the turning points do not occur at the same times.

Again, there is a connection here between the maximum steepness of the real solution and the trajectory of closest pole.  As has been noted, this may not be at all surprising, but it is interesting to see the numerical comparison, which suggests that, at least for our problem, a marker for whether the maximum steepness of the real solution is increasing or decreasing is to extend the solution to the complex plane and observe whether the distance of the closest singularity to real axis is increasing or decreasing.

\section{Large-time analysis}\label{sec:largetime}

\subsection{Numerical solution}

\begin{figure}[h!]
\centering
\includegraphics[width=\textwidth]{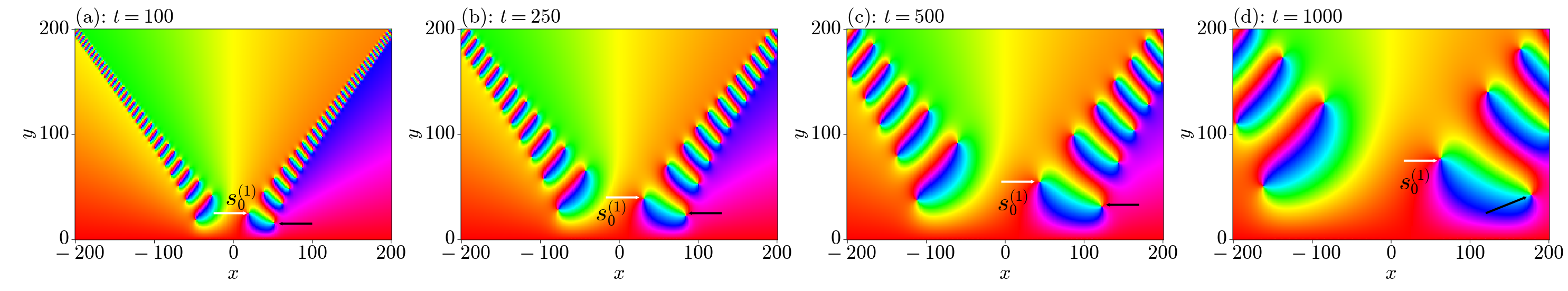}
\caption{Phase portraits of $u(z,t)$ for $\mu=1$, computed for times $t=100$, $250$, $500$ and $1000$. The closest pole to the real axis, $z=s^{(1)}_0$, and the associated zero are indicated by the white and black arrows, respectively.}

\label{fig:largishtime}
\end{figure}

Before deriving an asymptotic solution in the large-time limit, it is worth noting some key features of the solution for large times. In figure~\ref{fig:largishtime}, phase portraits and analytical landscapes of the solution are shown for $t=100$, $250$, $500$ and $1000$.  The poles appear to line up very closely to the rays $\arg(z)=\pi/4$, $3\pi/4$, and also appear to be moving further apart as $t$ increases.  As noted previously, for each pole there appears to be an associated simple zero.  Interestingly, these zeros appear to be propagating away from the poles so that the distance between each pole and zero in a pole-zero pairing is increasing in time.

\subsection{Similarity solution}

For large time, the leading-order analysis is standard, but we include some details here for completeness. Given the mass-preserving property of Burgers' equation (\ref{eq:maineq}), we look for a similarity solution of the form
\begin{equation}\label{eq:largetimeansatz}
u\sim \frac{\sqrt{\mu}}{t^{1/2}} \,\Psi(\eta)\quad \text{as}~t \to \infty,
\end{equation}
where
$$
\eta=\frac{x}{\sqrt{\mu t}}.
$$
Substituting this ansatz \eqref{eq:largetimeansatz} into Burgers' equation (\ref{eq:maineq}) leads to the second-order ordinary differential equation (ode)
\begin{equation}
\label{eq:intermediateansatzlargeimte}
- \sfrac{1}{2}\Psi -\sfrac{1}{2}\eta\Psi' + \Psi\Psi' = \Psi'',
\end{equation}
where the primes mean derivatives with respect to $\eta$.  Note this equation is the same as (\ref{eq:innerpdeconstantterm}), except that now we have scaled $\mu$ out.  The key difference here relates to the boundary conditions.  We integrate (\ref{eq:intermediateansatzlargeimte}) directly, noting that the resulting integration constant is zero if we enforce the condition that $\Psi$, $\Psi^2$, and $\Psi'$ all decay exponentially as $x \to \pm \infty$ (in contrast, after we integrate (\ref{eq:innerpdeconstantterm}), we end up with a nonzero constant by matching with (\ref{eq:farfield2})). Hence, we arrive at the Riccati equation
\begin{equation}
\label{eq:ricattilargeimteansatz}
\Psi^2-\eta\Psi=2\Psi',
\end{equation}
whose well-known exact solution is
\begin{equation}\label{eq:psiformlargetime}
\Psi = \frac{2\e^{-\eta^2/4}}{\sqrt{\pi}(\gamma - \erf\left({\eta}/{2}\right))}
\end{equation}
for some constant $\gamma>1$.

The function in (\ref{eq:psiformlargetime}) represents a one-parameter family of solutions that depend on $\gamma$.  To fix $\gamma$, we can use conservation of mass, noting that if we let $M = \int_{-\infty}^\infty u(x, t) \dx$, then $\mathrm{d}M/\mathrm{d}t=0$.  We write
\begin{align*}
M &  \sim \mu\int_{-\infty}^\infty \Psi(\eta) \dd{\eta}=2\mu\log\left(\frac{\gamma+1}{\gamma-1}\right)
\quad \text{as}~t\to\infty.
\end{align*}
For a given initial condition $u(x,0)=u_0(x)$, we also have $M = \int_{-\infty}^\infty u(x,0) \dx$, which provides a relationship between $\gamma$ and the initial profile.  For our case (\ref{eq:ic}), $M= \pi$, which means that
\begin{equation}
\gamma = \coth\left(\frac{\pi}{4\mu}\right).
\label{eq:gamma}
\end{equation}
Our similarity solution is therefore (\ref{eq:psiformlargetime}) with (\ref{eq:gamma}).
In the limit $\mu\rightarrow 0$, $\gamma\rightarrow 1^+$ and (\ref{eq:psiformlargetime}) reduces to the small-viscosity approximation given in the appendix of Chapman et al.~\textcite{Chapman_2007}.  Note that a different initial condition for our original problem with the same value of $M$ will lead to the same leading-order approximation in this large-time limit.

\begin{figure}[h!]
\centering
\includegraphics[width=\textwidth]{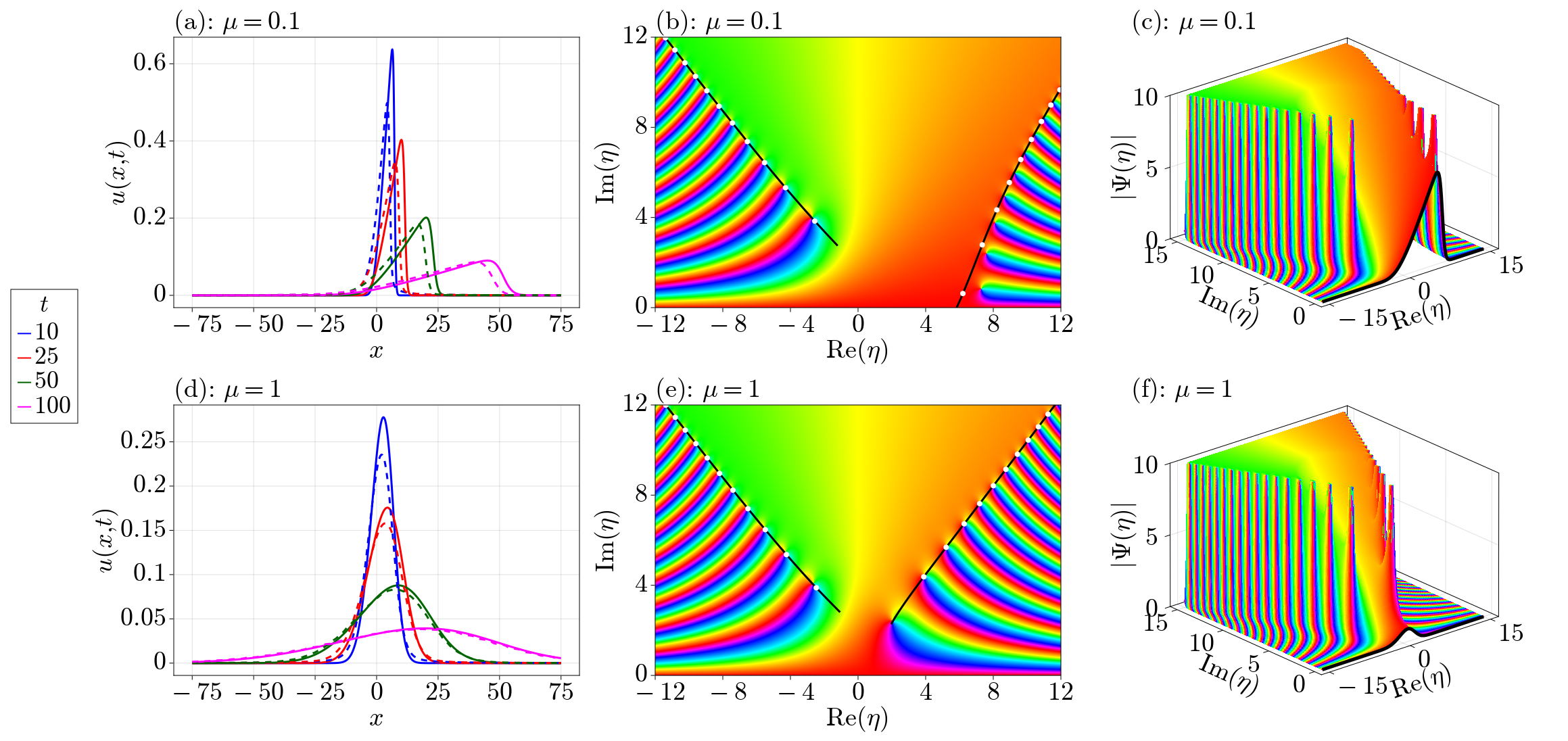}
\caption{Similarity solution for large time.  (a) and (d) show a comparison of  \eqref{eq:largetimelimit} (solid line) with the corresponding exact solution (dashed line) at times $t = 10, 25, 50, 100$ for $\mu=0.1$ and $1$, respectively.
(b) and (e) show the phase portraits of \eqref{eq:largetimelimit} in the $\eta$-plane, while (c) and (f) show the analytical landscapes.  The black curves in (b) and (e) are drawn using the approximations (\ref{eq:thetarhofarfield}) and (\ref{eq:thetarhofarfield2}), while the solid (white) dots are found by solving the transcendental equations (\ref{eq:firstquadtrans}) and (\ref{eq:secondquadtrans}).
}\label{fig:largetimerealcomps}
\end{figure}

In terms of the original variables, this large-time approximation in the complex plane becomes
\begin{equation}
\label{eq:largetimelimit}
u\sim
\frac{2\sqrt{\mu}\,\e^{-z^2/4\mu t}}{\sqrt{\pi t}
\left(\coth\left({\pi}/{4\mu}\right) - \erf\left({z}/{2\sqrt{\mu t}}\right)\right)}
\quad \text{as}~t\to\infty.
\end{equation}
In figure~\ref{fig:largetimerealcomps}(a) and (d) we compare exact solutions with (\ref{eq:largetimelimit}) on the real line for $\mu=0.1$ and $1$.  We again see that for small $\mu$ the solutions are steep at the front and for larger $\mu$ the solutions more closely resemble a Gaussian curve.
We also show the phase portraits and analytical landscapes for $\Psi(\eta)$ from (\ref{eq:psiformlargetime}) with (\ref{eq:gamma}), again for $\mu=0.1$ and $1$, again showing only the upper-half plane.  We observe the array of poles in the first and second quadrants that tend to lie on the rays $\arg(\Psi)=\pi/4$, $3\pi/4$ in the far field.  Further, we note each pole no longer has an associated zero in the plane, but instead the zeros in the first quadrant appear to approach $x=+\infty$ as $t\rightarrow\infty$, while the zeros in the second quadrant appear to approach $x=-\infty$.

\subsection{Locating poles of $\Psi(\eta)$ in the far field}

While the large-time analysis of Burgers' equation could be viewed as an undergraduate exercise, any analysis of the resulting singularities is much less common.  As such, we use the similarity solution to explore the geometry of poles of our solution for $u$ for late times.  Taking the leading-order expression \eqref{eq:largetimelimit}, in the upper-half plane we have $\erfc(\eta/2)\sim C+(2/\sqrt{\pi}\eta)\,\e^{-\eta^2/4}$ as $\eta\rightarrow\infty$, where one way to interpret this asymptotically is to set $C=0$ in the first quadrant and $C=2$ in the second quadrant.  Therefore the far-field behaviour of (\ref{eq:psiformlargetime}) is
\begin{equation}\label{eq:xifarfield}
\Psi(\eta) \sim \left\{
\begin{array}{ll}
\displaystyle \frac{2\,\e^{-\eta^2/4}}{\sqrt{\pi}(\gamma-1)}, & 0\leq\arg(\eta)\leq\pi/4 \\
\displaystyle \eta, & \pi/4<\arg(\eta)<3\pi/4 \\
\displaystyle \frac{2\,\e^{-\eta^2/4}}{\sqrt{\pi}(\gamma+1)}, & 3\pi/4\leq\arg(\eta)\leq\pi,
\end{array}
\right.
\end{equation}
which means that the poles must lie very slightly off the rays $\arg(\eta)=\pi/4$ and $3\pi/4$ and align themselves increasingly closer to these rays as $|\eta|$ increases.

To locate these poles we shall look for zeros of the denominator
\begin{equation}\label{eq:largetimeroot}
G(\eta) = \gamma - \erf\left(\frac{\eta}{2}\right) = \gamma - 1 + \erfc\left(\frac{\eta}{2}\right),
\end{equation}
where $\gamma$ is given by (\ref{eq:gamma}). By writing $\eta=\rho\,\e^{\i\theta}$, close to the ray $\arg(\eta)=\pi/4$ we have
$$
G\sim \gamma - 1 + \frac{2}{\sqrt{\pi}\rho}\exp\left(-\sfrac{1}{4}\rho^2\cos 2 \theta\right)\left[\cos\left(\sfrac{1}{4}\rho^2\sin 2 \theta + \theta\right) - \i \sin\left(\sfrac{1}{4}\rho^2\sin 2 \theta + \theta\right)\right]
\quad
\mbox{as}\quad \rho\rightarrow\infty.
$$
For the imaginary part of $G$ to vanish we require
\begin{equation}
\sin\left(\sfrac{1}{4}\rho^2\sin 2 \theta + \theta\right)=o(1),
\quad
\cos\left(\sfrac{1}{4}\rho^2\sin 2 \theta + \theta\right)\sim -1,
\quad
\mbox{as}\quad \rho\rightarrow\infty.
\label{eq:sincosfareta}
\end{equation}
By expanding these about $\theta=\pi/4$, we find $\rho^2/4+\pi/4\sim (2n+1)\pi$, where $n$ is an integer.
Using the results in (\ref{eq:sincosfareta}), for the real part to vanish we require
$$
\gamma - 1 - \frac{2}{\sqrt{\pi}\rho}\exp\left(-\sfrac{1}{4}\rho^2\cos 2 \theta\right)=o(1),
$$
or
\begin{equation}
\theta\sim \frac{\pi}{4}+ \frac{2}{\rho^2}\left(\ln\rho+\ln\left(\sfrac{1}{2}\sqrt{\pi}(\gamma-1)\right)\right)
\quad
\mbox{as}\quad \rho\rightarrow\infty.
\label{eq:thetarhofarfield}
\end{equation}
The poles of $\Psi$ in the first quadrant lie approximately on this path, which is indicated by the black curves in figure~\ref{fig:largetimerealcomps}(b) and (e).  Further, by combining with (\ref{eq:sincosfareta}), we find that values of $\rho$ to use in (\ref{eq:thetarhofarfield}) are found by solving the transcendental equation
\begin{equation}
-\sfrac{1}{2}\rho^2\tan\left(\sfrac{1}{4}\rho^2+\sfrac{\pi}{4}\right)
=\ln\rho+\ln\left(\sfrac{1}{2}\sqrt{\pi}(\gamma-1)\right).
\label{eq:firstquadtrans}
\end{equation}
We can further approximate solutions to this equation by
\begin{equation}
\rho^2 \sim \left(8n+3\right)\pi-\frac{\ln 2n + 2\ln(\pi(\gamma-1))}{2n\pi}
\quad\text{as}~n\to\infty.\label{eq:firstquadrho}
\end{equation}
Thus we can provide good evidence there is a countably infinite number of these poles in the first quadrant and give asymptotic descriptions for their locations, together with their separation distance which decreases like $\sqrt{2\pi/n}$ in $\eta$-plane as $n\rightarrow\infty$ (which is $\sqrt{2\mu\pi/n}\,t^{1/2}$ in the $z$-plane).  We use (\ref{eq:firstquadrho}) to draw the solid white dots in the first quadrant of figure~\ref{fig:largetimerealcomps}(b) and (e), which appear to line up with the poles very well; note the predictions for $\mu = 0.1$ take longer to settle onto the pole locations than in the $\mu = 1$ case, but numerical tests show they do indeed provide very accurate approximations as $|\eta|$ increases beyond the scale of this figure.

An almost identical analysis applies in the second quadrant.  {
The expansion for $G$ close to the ray $\arg(\eta)=3\pi/4$ is the same as above, mutatis mutandis; however, again we arrive at equation (\ref{eq:sincosfareta}).  This time, expanding (\ref{eq:sincosfareta}) about $\theta=3\pi/4$, we find $\rho^2/4+\pi/4\sim 2n\pi$, where $n$ is an integer.}
Leaving out further details, we find all the poles of $\Psi$ lie approximately on the path
\begin{equation}
\theta\sim \frac{3\pi}{4}- \frac{2}{\rho^2}\left(\ln\rho+\ln\left(\sfrac{1}{2}\sqrt{\pi}(\gamma+1)\right)\right)
\quad
\mbox{as}\quad \rho\rightarrow\infty,
\label{eq:thetarhofarfield2}
\end{equation}
which is drawn on figure~\ref{fig:largetimerealcomps}(b) and (e) as a black curve.  The values of $\rho$ along this curve are determined as solutions to the transcendental equation
\begin{equation}
-\sfrac{1}{2}\rho^2\tan\left(\sfrac{1}{4}\rho^2+\sfrac{\pi}{4}\right)
=\ln\rho+\ln(\sfrac{1}{2}\sqrt{\pi}(\gamma+1)),
\label{eq:secondquadtrans}
\end{equation}
which leads to
\begin{equation}
\rho^2 \sim \left(8n-1\right)\pi-\frac{\ln 2n + 2\ln(\pi(\gamma+1))}{2n\pi}
\quad\text{as}~n\to\infty.\label{eq:secondquadrho}
\end{equation}
Again, these predictions, indicated by white dots in the second quadrant in figure~\ref{fig:largetimerealcomps}(b) and (e), line up extremely well with the poles.  The separation distance in the second quadrant again decreases like $\sqrt{2\pi/n}$ as $n\rightarrow\infty$.

\subsection{Anti-Stokes line analysis}

We now complement the above results by providing a more general, but less detailed, description of the far-field structure of the singularities applicable for $t=\OO(1)$ or $t\gg 1$ (the case $t\rightarrow 0^+$ with $|z|=\OO(1)$ is likely to be amenable to similar methods, but will be rather more involved).  Little of what follows depends on the specific initial conditions in question.

There are two contributions to the $z\rightarrow +\infty$ limit, namely the divergent algebraic series and the exponential term that appears beyond all algebraic orders.  The former is determined by substituting the ansatz
$$
u_A \sim \sum_{n=0}^\infty \frac{b_n(t)}{z^n}
\quad\mbox{as}\quad z\rightarrow\infty
$$
into (\ref{eq:maineq}) and matching back onto
$$
u_A \sim \frac{1}{z^2}-\frac{1}{z^4}+\frac{1}{z^6}+\ldots
$$
(the initial condition (\ref{eq:ic})) as $t\rightarrow 0^+$.  The result is
\begin{equation}
u_A\sim \frac{1}{z^2}+\frac{6\mu t-1}{z^4}+\frac{2t}{z^5}+\ldots.
\label{eq:jrkalgebraic}
\end{equation}
The exponential term results from this algebraic series turning on the Liouville-Green (WKB) contribution
\begin{equation}
u_{LG}\sim \frac{1}{t^{1/2}}\phi\left(\frac{z}{t}\right)\e^{-z^2/4\mu t}
\label{eq:jrkLG}
\end{equation}
across the Stokes line $\Im(z)=0$ (this can be interpreted as corresponding to the merging in the far field of the two Stokes lines $\Im(z)=\pm 1$ that arise in the small-time analysis, as suggested in section~\ref{sec:smalltimeWKB}); in (\ref{eq:jrkLG}), $\phi(\zeta)$ is an arbitrary function whose large-$\zeta$ and small-$\zeta$ behaviour can be determined from the small- and large-time analyses, respectively, but whose detailed form has little bearing on what follows.

The anti-Stokes line in $\Im(z)>0$, on which (\ref{eq:jrkalgebraic}) and (\ref{eq:jrkLG}) become comparable, has $\arg(z)=\pi/4$ and we set
$$
z=r\e^{\i\pi/4}+\frac{\sigma\,\e^{3\i\pi/4}}{r}
$$
with $r\gg 1$ fixed and $\sigma\in\mathbb R$ spanning the anti-Stokes line (refinements to this scaling are required below).  Thus
\begin{equation}
u_A\sim -\frac{\i}{r},
\quad
u_{LG}\sim \frac{1}{t^{1/2}}\phi\left(\frac{r\e^{\i\pi/4}}{t}\right)\e^{-\i r^2/4\mu t}
\,\e^{\sigma/2\mu t}.
\label{eq:jrkALG}
\end{equation}
To make these two expressions comparable for $r\gg 1$, we translate and rescale $\sigma$ in the form
$$
\sigma = \mu t\ln t-2\mu t\ln \phi(r\e^{\i\pi/4}/t)-4\mu t\ln r+2\mu t\bar{\sigma},
$$
so that
\begin{equation}
u_{LG}\sim \frac{1}{r^2}\e^{-\i r^2/4\mu t}
\,\e^{\bar\sigma}.
\label{eq:jrknewLG}
\end{equation}
Now guided by (\ref{eq:jrknewLG}) we introduce the new time variable $\tau=-r^2/4\mu t$.  Retaining the dominant contributions to each of the three terms in the pde renders
$$
\frac{\partial u}{\partial\tau}+\frac{2t}{r\e^{3\i\pi/4}}\,
u\,\frac{\partial u}{\partial\bar\sigma}\sim \i\,\frac{\partial^2 u}{\partial\bar\sigma^2},
$$
so for $u=\OO(1/r^2)$, as required by (\ref{eq:jrkALG}) and (\ref{eq:jrknewLG}), the nonlinear term is negligible and
\begin{equation}
u\sim \frac{1}{r^2}(-\i+\e^{\bar\sigma+\i\tau}),
\label{eq:jrkorderu2}
\end{equation}
which matches with (\ref{eq:jrkALG}) and (\ref{eq:jrknewLG}).  That is, due to linearity at leading order, the Liouville-Green contribution simply passes through the algebraic series, coming to dominate as $\sigma\rightarrow\infty$.

In order to bring the convective term into play, the further translation
$$
\bar\sigma = 3\ln r-\ln t+\hat\sigma
$$
is required and setting $u=rv/t$ then implies the leading-order balance,
$$
\frac{\partial v}{\partial\tau}-2\,\e^{\i\pi/4}
v\,\frac{\partial v}{\partial\hat\sigma} = \i\,\frac{\partial^2 v}{\partial\hat\sigma^2},
$$
leading to the travelling-wave solution
\begin{equation}
v=\frac{1}{e^{-(\hat\sigma+\i\tau)}+\e^{\i\pi/4}}
\label{eq:jrktw}
\end{equation}
on matching with (\ref{eq:jrkorderu2}).  Asymptotic expressions for the singularity locations as $r\rightarrow\infty$ can thus be determined as $\sigma=0$ with
\begin{equation}
\frac{r^2}{4\mu t}\sim \left(2n+\sfrac{3}{4}\right)\pi
\label{eq:jrkpolelocation}
\end{equation}
for positive integer $n\gg 1$, which agrees with (\ref{eq:firstquadrho}).  The result (\ref{eq:jrkpolelocation}) explicitly relates the instantaneous generation of an infinite number of singularities present for arbitrary large $r$ to the infinite-speed-of-propagation property of the heat operator, and confirms that they are more-and-more closely spaced as $r\rightarrow\infty$.

Note that a similar analysis can be applied in the second quadrant to show the singularity locations behave like $r^2/4\mu t\sim \left(2n-\sfrac{1}{4}\right)\pi$ as $r\rightarrow\infty$ (as with (\ref{eq:secondquadrho})), but we do not include the details here.

\section{Rational approximations via the AAA algorithm}\label{sec:aaa}

Clearly we have made much progress due to the exact solution to Burgers' equation.  Without such an exact formula, we would be left with asymptotics and other analytical approximations.  Here we briefly explore the use of a rational approximation that is computed numerically using the AAA (adaptive Antoulas--Anderson) algorithm \autocite{nakatsukasa2018aaa}.  We will show how we can view an analytic continuation of a solution in the complex plane from numerical data on the real line, and then we will show how we can track the closest pole to the real axis, just as we did in section~\ref{sec:poletrajectory}.

The AAA algorithm takes in some points $\{z_j\}$ and corresponding function values $f_j$, and returns a rational function $r(z)$ that can be used to approximate the original function for values $z \in \mathbb C$. In particular, we can use this algorithm to provide an approximate analytic continuation of solutions of \eqref{eq:maineq} from a numerical solution on the real line, allowing us to view the function and its poles in the complex plane. We note that (numerical) analytic continuation is ill-posed and becomes much less accurate the further we are from the real line \textcite{trefethen2019approximation, trefethen2020quantifying}. We solve the pde \eqref{eq:maineq} using the \texttt{pde15s} function in Chebfun \autocite{driscoll2014chebfun}, and the AAA algorithm is implemented in Chebfun's \texttt{aaa} function.

\begin{figure}[h!]
\centering
\includegraphics[width=\textwidth]{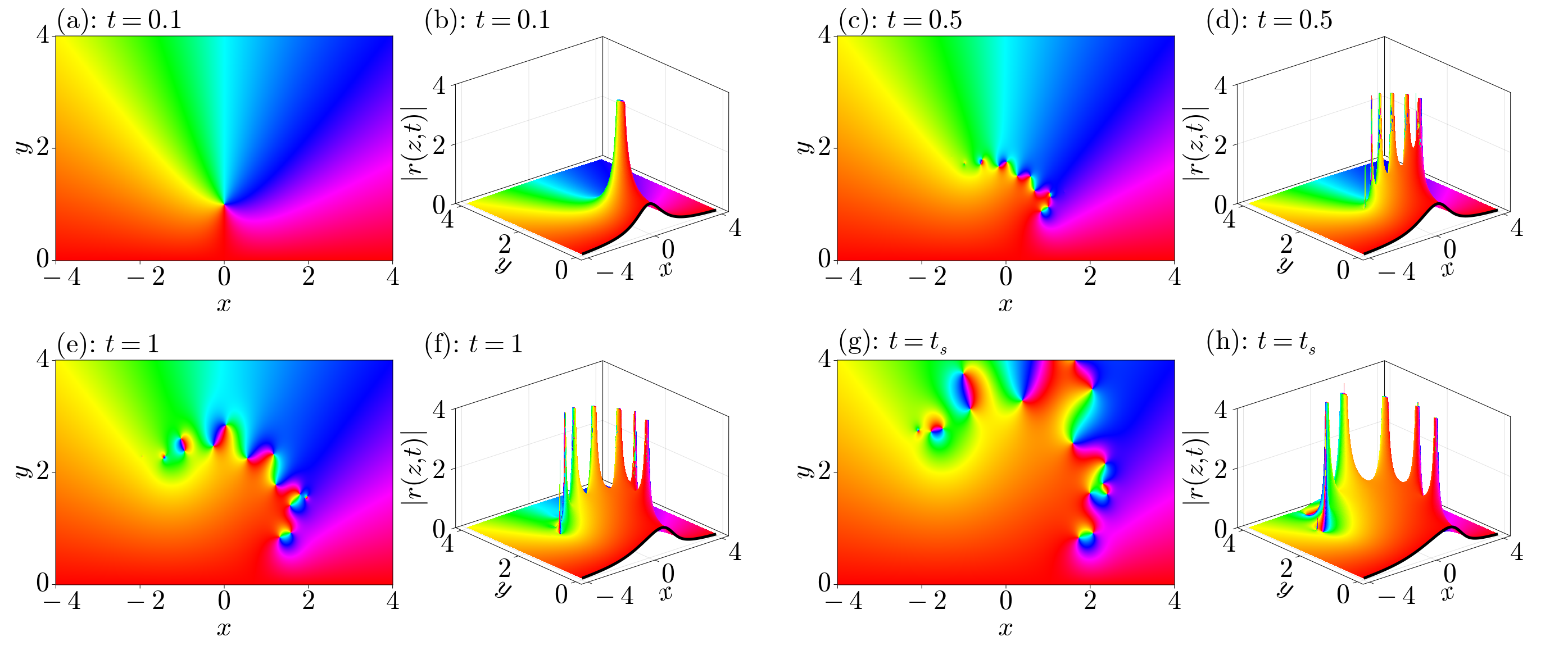}
\caption{Phase portraits and analytical landscapes of the AAA approximant $r(z, t)$ using data from a numerical solution of \eqref{eq:maineq} on the real line at $\mu = 0.1$.}\label{fig:aaaviscousplots1}
\end{figure}

Figure \ref{fig:aaaviscousplots1} shows plots of the AAA approximant $r(z, t)$ that approximate $u(z, t)$ in the complex plane for the case $\mu=0.1$.  This approximant is formed using data on the real line on the domain $-15 \leq x \leq 15$ with $N =250$ spatial grid points. Comparing figure \ref{fig:aaaviscousplots1} to the corresponding plots in figure \ref{fig:complexsolutionplotviewsmallmu}, we see that the AAA algorithm does not produce the same array of poles in the first and second quadrants as we move far away from the real axis, although predicting far-field behaviour like this is known to be very difficult to capture with numerical analytic continuation \autocite{trefethen2020quantifying}.

\begin{figure}[h!]
\centering
\includegraphics[width=0.8\textwidth]{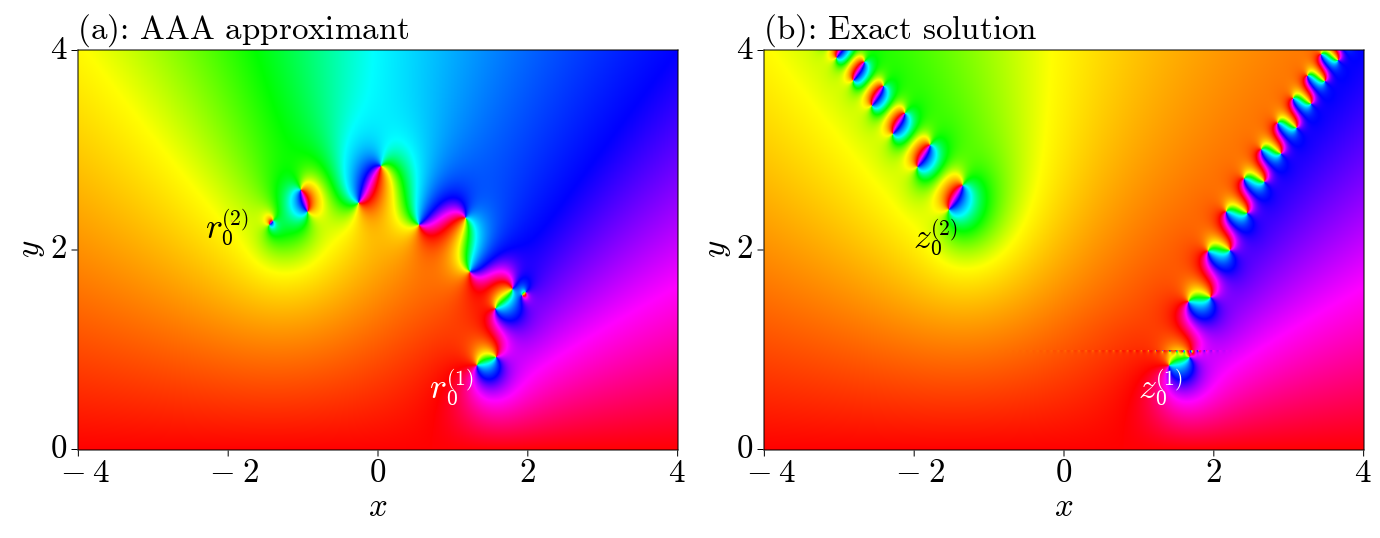}
\caption{Plots of the (a) AAA approximant $r(z, t)$ and the (b) exact solution at $\mu=0.1$ and $t=1$. }
\label{fig:aaaviscousplots1specific}
\end{figure}

On the other hand, there are key properties of the solution that AAA picks up very well.  For example, the colouring in figures \ref{fig:aaaviscousplots1} and figure \ref{fig:complexsolutionplotviewsmallmu} is very similar in a strip around the real line and further into the complex plane for sufficiently large $\Re(z)$, suggesting that AAA provides a good prediction of phase on this part of the domain.  Comparing the solutions further, we see that AAA appears to have picked up the first two poles in the first quadrant and the first pole in the second quadrant quite well.  To make this observation clearer, in figure~\ref{fig:aaaviscousplots1specific} we show the $t=1$ solution from the AAA approximant in (a) and the exact solution in (b). In these plots, we find that $s^{(1)}_0$ and $s^{(1)}_1$ are very well approximated by $r^{(1)}_0$ and $r^{(1)}_1$ in AAA, while $s^{(2)}_0$ is reasonably well approximated by $r^{(2)}_0$.

We now return to figure \ref{fig:saddlepointtrajectorytracked}, where we plotted trajectories of the pole in the first quadrant closest to the real line, namely $z=s^{(1)}_0$, for three different values of viscosity using the exact solution (and also an approximation from steepest descents).  Also included in that figure are data computed using the AAA algorithm, which we see matches the other two curves extremely well. The curves here are obtained using Chebfun's \texttt{pde15s} and \texttt{aaa} functions \autocite{driscoll2014chebfun}, and the algorithm we use for tracking these poles is described in more detail in \ref{sec:appendixaaa}. We see that, for all $\mu$ values shown, we can reliably track this pole. Thus, in our problem the AAA algorithm provides a reliable numerical method for tracking the closest pole to the real axis, and presumably would be equally useful for a host of other nonlinear pdes (including those that do not have exact solutions).

\section{Discussion}\label{sec:discussion}

In this work, we present an analysis of solutions of Burgers' equation \eqref{eq:maineq} in the complex plane.  A nice feature of Burgers' equation is that all singularities must be simple poles, leading to a simpler analysis than if we had to consider other singularities such as branch points. While previous work has considered initial conditions that are entire and are odd on the real line, leading to poles that are constrained to the imaginary axis and coming in from infinity  \cite{bessis1984pole,bessis1990complex,senouf1997dynamicsA, senouf1997dynamics,caflisch2015complex,weideman2021dynamics}, we use the initial condition (\ref{eq:ic}) which is neither odd on the real line nor entire in the plane. This initial condition therefore allows us to study the two-dimensional trajectories of poles in the complex plane as they emerge from $z = \pm \i$.  Note that Chapman et al.~\cite{Chapman_2007} provide a comprehensive analysis of \eqref{eq:maineq} with (\ref{eq:ic}) in limit $\mu\rightarrow 0^+$; our study is different in that we are interested mostly in $\mu = \OO(1)$ and therefore our work can be seen as complementary to that in~\cite{Chapman_2007}.

Our analysis provides a clear picture of the dynamics of the complex singularities of  \eqref{eq:maineq} with \eqref{eq:ic}. We hypothesise that there is a countably infinite number of simple poles that emerge spontaneously from $z=\i$ at $t=0$ and align themselves in an array in both the first and second quadrants of the complex plane.  The form of a small-time expansion around $z = \i$ leads to the inner problem  \eqref{eq:innerpde}  with an exact solution \eqref{eq:phisolutionparaboliccylinder} in terms of parabolic cylinder functions. We find that the closest pole to the real line, $z=s_0^{(1)}(t)$, initially moves as $s_0^{(1)}(t) \sim \i + t^{1/2}\xi_0$ as $t \to 0^+$, where $\xi_0$ is the pole of \eqref{eq:phisolutionparaboliccylinder} with the lowest imaginary part.  For sufficiently small $\mu$, $\xi_0$ lies in the fourth quadrant of the $\xi$-plane, which means that $s_0^{(1)}(t)$ initially heads towards the real $z$-axis.  In contrast, for sufficiently large $\mu$, all poles of $u$ move away from the real axis.  Finally, by considering the far-field behaviour of the inner solution \eqref{eq:phisolutionparaboliccylinder}, we derive the asymptotics \eqref{eq:Phifarfield} with different behaviours in different sectors of the complex $\xi$-plane (due to the Stokes phenomenon).  This result  leads to the interesting observation that the analytic continuation of the solution of \eqref{eq:maineq} with \eqref{eq:ic} does not approach $u(z, 0) = 1/(1+z^2)$ everywhere in the $z$-plane as $t \to 0^+$.

We are able to track the closest pole to the real line for small to moderate times reliably by applying rootfinding to the Cole-Hopf solution or the method of steepest descents, or from implementing a rational approximation with the AAA algorithm.  Our analysis allows us to make connections between the trajectory of the closest pole to the real line and the maximum steepness of the solution profile on the real line. For example, we observe that the solution on the real line initially steepens if the closest pole to the real line initially moves closer to the real line or initially flattens out when the pole initially moves away.  (These options for small-time behaviour are not present for initial conditions that are entire  \cite{bessis1984pole,bessis1990complex,senouf1997dynamicsA, senouf1997dynamics,caflisch2015complex,weideman2021dynamics}, since in that situation poles always initially move closer to the real axis, given they come in from $|z|=\infty$.)  These connections are not surprising, and are reminiscent of links made between oscillations on the real line and the closest singularity for dispersive systems \cite{weideman2021dynamics,senouf1996}.

Another observation is that for each simple pole there appears to be an associated simple zero, which is very close to the pole for small time but moves further away from it for large time.  By deriving a similarity solution of the form \eqref{eq:largetimeansatz} as $t \to \infty$, we are able to determine the leading-order behaviour of the solution for large time \eqref{eq:largetimelimit}.  We find that for large time, the poles in the first and second quadrant must lie very close the rays at angles $\pi/4$ and $3\pi/4$.  Moreover, we find that the spacing between these poles in both the first and second quadrants is $\sqrt{2\mu \pi t/n}$ as $n \to \infty$, which is also the scaling for $t\ll 1$ and $t=\OO(1)$.  Finally, our large-time limit suggests each zero in a pole-zero pairing evolves towards $z = \pm \infty$ as $t \to \infty$.

Much of the analysis in this work relies on the knowledge of exact solutions, including the Cole-Hopf solution to the full problem (\ref{eq:exactsolutionrunge}), the solution to the small-time inner problem in terms of parabolic cylinder functions (\ref{eq:phisolutionparaboliccylinder}), or the solution to the large-time leading-order problem in terms of the error function (\ref{eq:psiformlargetime}).  Without these exact solutions, we have successfully employed the AAA algorithm for tracking the closest pole to the real axis for $\mu = \OO(1)$ using only a numerical solution on the real line.  We have also provided a variety of asymptotic results that do not rely on exact solutions.  And further, as mentioned above, significant analytical progress can be made without exact solutions via the limit $\mu\rightarrow 0^+$ \cite{Chapman_2007}. It would be interesting to apply asymptotic methods in other parts of our study to extract meaningful information, again without knowledge of the exact solution.  For example, we could go further than the Liouville-Green analysis in \ref{eq:appendixinner} and apply techniques in exponential asymptotics to study how the exponentially small term $C\xi^{-\i/2\mu}\e^{-\xi^2/4\mu}$ is switched on across the Stokes lines (the rays $\arg(\xi)=0$, $\pi$) for the small-time inner problem (\ref{eq:innerpdeconstantterm}) with (\ref{eq:farfield2}).  Further analysis to approximate the location of the poles in the $\xi$-plane could be conducted by resummation of the appropriate transseries, cf.~\cite{costin2001}.  In our large-time limit, similar calculations could be attempted without the exact solution.  We leave these issues for further study.

Finally, we emphasise that the small-time analysis we have conducted with the initial condition (\ref{eq:ic}), which has simple poles, is much cleaner than the cases for which the initial condition has other types of singularities. For example, for the initial condition
\begin{equation}\label{eq:ic2}
u(x, 0) = u_0(x) = \frac{1}{(1+x^2)^\beta},
\end{equation}
where $\beta>0$, there are three distinct regimes.  As just mentioned, $\beta=1$ is the borderline case considered in section~\ref{sec:smalltime}.  Notably, with $\beta=1$, the advective term $u_0u'_0$ balances the diffusive term $\mu u''_0$ near $x=\mathrm{i}$, which leads to a leading-order problem (\ref{eq:innerpdeconstantterm}) whose solutions have simple poles (albeit in general with different magnitude than those of the initial data). For $0<\beta<1$, the initial condition has a branch point at $x=\mathrm{i}$.  Here, the diffusive term dominates, and so to leading order the small-time behaviour appears to be governed by a second-order linear ode whose solutions do not have singularities.  A further rescaling is required in order to balance the nonlinearity, leading to infinitely many simple poles, as expected.  We summarise these ideas in \ref{eq:gammanot1}.  For $\beta>1$, the initial condition has either a branch point or higher order poles (when $\beta>1$ is an integer) at $x=\mathrm{i}$.  In this case, the advective term dominates the diffusive term near $x=\mathrm{i}$, leading to a first-order nonlinear problem with square-root branch points.  A further rescaling is needed to link in with the higher-order diffusion term, which gives rise to a second-order nonlinear ode whose solutions have infinitely many simple poles.  Again, we summarise the main ideas for $\beta>1$ in \ref{eq:gammanot1}.  All of this analysis for $\beta\neq 1$ is much more complicated than that for $\beta=1$; however, we emphasise that for $t=\mathcal{O}(1)$, the qualitative behaviour of the solutions, with arrays of poles in each quadrant, is the same regardless of $\beta$.

\section*{Acknowledgements}
DJV is grateful for the financial support of a QUT Vice-Chancellor's Scholarship (Academic) in 2021, and from QUT's Centre for Data Science for a QUT Centre for Data Science Scholarship in 2022.
JRK gratefully acknowledges a Fellowship from the Leverhulme Trust.  SWM, CJL and JRK would like to acknowledge the Isaac Newton Institute for Mathematical Sciences, Cambridge, for support during the programme Applicable Resurgent Asymptotics: Towards a Universal Theory where part of the work on this paper was undertaken.  This programme was supported by the EPSRC grant EP/R014604/1.



  \bibliographystyle{elsarticle-num}
  \bibliography{VandenHeuvel2022_Burgers_Physica_References.bib}

%
%

\appendix

\section{Small-time analysis with the saddle-point method}\label{sec:smalltimesaddle}

To analyse \eqref{eq:maineq} in the limit $t \to 0^+$ using the method of steepest descents, we start by writing \eqref{eq:exactsolutionrunge} in the form
\begin{equation}
u(z, t) = \frac{N(z, t)}{D(z, t)} = \dfrac{\displaystyle\int_{-\infty}^\infty \frac{z-s}{t}G(s)\e^{g(s)/t} \ds}{\displaystyle\int_{-\infty}^\infty  G(s)\e^{g(s)/t} \ds},
\quad
G(s) = \exp\left\{-\frac{1}{2\mu}\arctan(s)\right\},
\quad g(s)=-\frac{(z-s)^2}{4\mu},
\end{equation}
where
\begin{equation}
 G(s) = \exp\left\{-\frac{1}{2\mu}\arctan(s)\right\}\quad \text{and} \quad g(s)=-\frac{(z-s)^2}{4\mu},\label{eq:denominatorintegralsaddlepoint}
\end{equation}
and then identify that the single saddle point of $g(s)$ in \eqref{eq:denominatorintegralsaddlepoint} is at $s_0=z$. We deform onto this contour by writing $g(s)-g(s_0)=-\zeta^2$, so that $s=z+2\xi\sqrt{\mu}$ for $-\infty<\zeta<\infty$; then we obtain
$$
N(z, t) \sim -4\mu t^{1/2}\sum_{n=0}^\infty e_{2n+1}\Gamma\left(n+\frac32\right)t^n, \quad
D(z, t) \sim 2\mu^{1/2}t^{1/2}\sum_{n=0}^\infty e_{2n}\Gamma\left(n+\frac12\right)t^n,\quad\text{as}\quad t \to 0^+,
$$
where the $e_n$ are defined by $e_0 = \exp(-\arctan(z)/(2\mu))$ and
\begin{equation}
e_n = \frac{1}{n}\sum_{\ell=1}^n \ell d_\ell e_{n-\ell},
\quad
d_n = \frac{2^{n-1}\mu^{n/2-1}}{n}(-1)^n\sin^n\left(\frac{\pi}{2}-\arctan z\right)\sin\left[n\left(\frac{\pi}{2}-\arctan z\right)\right],
\quad n=1,2,\ldots.
\end{equation}
The derivation of this expansion is given in the Supplementary Material. Taking two terms in  the series for the numerator and the denominator, we can find for example that
\begin{equation}\label{eq:smalltimesaddlepoint}
u(z, t) \sim \frac{1+8\mu z}{1+z^2}- \frac{8\mu^2(4z^3+4z+t)}{4\mu z^4 + 8\mu z^2 + 4 \mu t z + 4 \mu + t}, \quad t \to 0^+.
\end{equation}
Expanding \eqref{eq:smalltimesaddlepoint} around $t = 0^+$ to get a polynomial expansion shows that this expansion is what we obtained directly from the pde \eqref{eq:naiveexpansion} in \eqref{eq:naiveexpansionterms}, provided $|\mathrm{Im}(z)| < 1$.  For $|\mathrm{Im}(z)| \geq 1$, we need to deform the contour around the branch point at either $z=\mathrm{i}$ or $-\mathrm{i}$, which will lead to additional considerations.

\section{Further details on the inner region analysis}\label{eq:appendixinner}

\subsection{Liouville-Green (WKB) analysis}

We start by considering the second-order ode (\ref{eq:innerpdeconstantterm}) and boundary condition (\ref{eq:farfield2}).  Linearising about this boundary condition, we write
$$
\Phi_0 \sim -\frac{\i}{2\xi}+F(\xi),
\quad \xi \to -\i\infty,
$$
so that $F$ satisfies the linear equation
\[
-\frac12F - \frac12\xi F' - \mu F'' = \frac{\i}{2\xi}F'-\frac{\i}{2\xi^2}F.
\]
Applying a Liouville-Green (WKB) ansatz $F\sim\,\mathrm{exp}(\phi_0+\phi_1+\ldots)$, where $|\phi_0|\gg |\phi_1|$ in the usual way, then
$$
\mu\phi_0^{'2}+\frac{1}{2}\xi\phi_0'=0,
\quad
2\mu\phi_0'\phi_1'+\mu\phi_0^{''}+\frac{1}{2}\xi\phi_1'+\frac{1}{2}=-\frac{\i}{2\xi}\phi_0'.
$$
The first of these equation has two possible solutions, either $\phi_0=\,\mathrm{constant}$ or $\phi_0=-\xi^2/4\mu+\,\mathrm{constant}$, which leads to $\phi_1=-\log\xi+\,\mathrm{constant}$ or $\phi_1=-(\i/2\mu)\log\xi+\,\mathrm{constant}$, respectively.

With these two linearly independent solutions for $F$, we find that
\[
F \sim K_1\,\frac{1}{\xi}+K_2\,\xi^{-\i/2\mu}\e^{-\xi^2/4\mu},
\quad \xi \to -\i\infty,
\]
for some constants $K_1$ and $K_2$.  Now $\e^{-\xi^2/4\mu}$ grows exponentially in this limit thus, in order to satisfy (\ref{eq:farfield2}), we must take $K_1=K_2=0$. Therefore we conclude that there are no degrees of freedom, and so the condition (\ref{eq:farfield2}) is acting as two boundary conditions, which is all that is required for our original second-order ode (\ref{eq:innerpdeconstantterm}).

\subsection{Exact solution of inner problem}

We treat the Riccati equation (\ref{eq:firstorderpde}) by introducing the new independent variable $g(\xi)$ via $\Phi_0 = -2\mu(g'/g)$, which leads to the linear ode
\begin{equation}\label{eq:appd146}
g''+\frac{\xi}{2\mu}g'-\frac{\i}{8\mu^2}g=0.
\end{equation}
Now use the change of variables $\zeta=-\xi^2/(4\mu)$ so that
\begin{equation}\label{eq:appd147}
\zeta\dv[2]{g}{\zeta}+\left(\frac12-\zeta\right)\dv{g}{\zeta} +\frac{\i}{8\mu}g=0.
\end{equation}
This is Kummer's equation with $b=1/2$ and $a=-\i/(8\mu)$ \autocite[Equation 13.2.1]{NIST:DLMF}. Thus, the solution to \eqref{eq:appd147} is
\begin{equation}\label{eq:kummerfunctionsolutions}
g(\zeta) = A\,M^*\left(-\frac{\i}{8\mu^2}, \frac12, \zeta\right) + B\,U^*\left(-\frac{\i}{8\mu^2}, \frac12, \zeta\right),
\end{equation}
where $A$ and $B$ are some constants and $M^*$ and $U^*$ are the Kummer functions defined in Equations 13.2.2 and 13.2.6 of \textcite{NIST:DLMF}, respectively. If we now use the asymptotic expansions for $M^*$ and $U^*$ (given in Equations 13.2.23 and 13.2.2 of \textcite{NIST:DLMF}) we can write
\begin{equation}\label{eq:kummerasymptoticexpansison}
M^*\left(-\frac{\i}{8\mu}, \frac12, -\frac{\xi^2}{4\mu}\right) \sim \frac{\sqrt{\pi}}{\Gamma\left(-\frac{\i}{8\mu}\right)}
\frac{2\sqrt{\mu}}{(-\xi^2)^{1/2}}
\left(-\frac{\xi^2}{8\mu}\right)^{-\i/\mu}\e^{-\xi^2/4\mu}, \quad U^*\left(-\frac{\i}{8\mu}, \frac12, -\frac{\xi^2}{4\mu}\right) \sim \left(-\frac{\xi^2}{4\mu}\right)^{\i/8\mu}, \quad \xi \to -\i\infty.
\end{equation}
These expansions in \eqref{eq:kummerasymptoticexpansison} can be differentiated to compute $\Phi_0 = -2\mu(g'/g)$, from which we find that in order to satisfy (\ref{eq:farfield2}), we choose $A=0$ and $B=1$. Using these values in \eqref{eq:kummerfunctionsolutions} and then using the relationship between $\Phi_0$ and $g$ gives
\begin{equation}\label{eq:phisolutionkummerfunctions}
\Phi_0 = \frac{\xi\left(\frac{\i}{8\mu}\right)U^*\left(1-\frac{\i}{8\mu}, \frac32, -\frac{\xi^2}{4\mu}\right)}{U^*\left(-\frac{\i}{8\mu}, \frac12, -\frac{\xi^2}{4\mu}\right)},
\end{equation}
where we used Equation 13.3.22 of \textcite{NIST:DLMF} to differentiate $U^*$. We can now use the relationship between the parabolic cylinder function $U$ and this Kummer function $U^*$ when the second argument of $U^*$ is $1/2$ or $3/2$ to give (\ref{eq:phisolutionparaboliccylinder}).

The parabolic cylinder function $U(a, z)$ is computed using the relationship \autocite{WeissteinParabolic}
\begin{align}
U(a, z) &= \cos\left[\pi\left(\sfrac{1}{4}+\sfrac{1}{2}a\right)\right]Y_1(a, z) - \sin\left[\pi\left(\sfrac{1}{4}+\sfrac{1}{2}a\right)\right]Y_2(a, z),
\end{align}
where
\begin{align}
Y_1(a, z) = \frac{1}{\sqrt{\pi}}\frac{\Gamma\left(\frac14-\frac12a\right)}{2^{a/2+1/4}}\e^{-z^2/4}
{}_1F_1\left(\sfrac{1}{2}a+\sfrac{1}{4};\sfrac{1}{2};\sfrac{1}{2}z^2\right),
\quad
Y_2(a, z) = \frac{1}{\sqrt{\pi}}\frac{\Gamma\left(\frac34-\frac12a\right)}{2^{a/2-1/4}}z\e^{-z^2/4}
{}_1F_1\left(\sfrac{1}{2}a+\sfrac{3}{4};\sfrac{3}{2};\sfrac{1}{2}z^2\right),
\end{align}
and the hypergeometric function ${}_1F_1(a;b;c)$ is computed using the \texttt{HypergeometricFunctions.jl} package in \textsc{Julia} \autocite{Slevinsky2021}.

\section{Tracking poles with the AAA algorithm}
\label{sec:appendixaaa}

In this section, we describe how we apply the AAA algorithm for tracking poles. The tools we use for this analysis are done in MATLAB with the Chebfun software \autocite{driscoll2014chebfun}. The first step in this analysis is to obtain a numerical solution for Burgers' equation. This is done by approximating \eqref{eq:maineq} with the problem
\begin{equation}\label{eq:aaaapproximatepde}
\pdv{u}{t} + u\pdv{u}{x} = \mu\pdv[2]{u}{x}, \quad -L \leq x \leq L,
\end{equation}
with boundary conditions $\partial u(\pm L, t)/\partial x = 0$ and initial condition $u(x, 0) = 1/(1+x^2)$. The choice of $L$ should be suitably large, and in this work we take $L = 15$ --- if other nonlinear pdes are considered, then depending on the decay rate of the solution as $|x| \to \infty$, $L$ may need to be larger.  This problem is easily solved using the Chebfun function \texttt{pde15s}, where we solve up to $t = T$, $T=2$, returning the solutions at times $t_j = (j-1)\Delta t$, $j=1,\ldots,M$, where $\Delta t = T/(M-1)$. We take $M=501$ in this work. Now, rather than returning the solution at equally spaced gridpoints, which could cause issues when interpolating the data with the AAA algorithm \cite{huybrechs2022AAA, platte2011impossibility} we will return the solution at each time $t_j$ at the points $x_i$, where
\begin{equation}
x_i = 5\cos\left(\frac{2i-1}{500}\pi\right),\quad i=1,2,\ldots,N=250.
\end{equation}
These points are the Chebyshev points on the interval $[-5, 5]$, and are a natural choice for the Chebyshev series returned from Chebfun's \texttt{pde15s} \autocite{driscoll2014chebfun}. We let $u_{ij}$ denote the approximate solution to \eqref{eq:aaaapproximatepde} at $x = x_i$ for $i=1,\ldots,N$, $t = t_j$ for $j=1,\ldots,M$, and for some $\mu$ value.

Now having the data from the numerical solution, we can obtain the AAA approximants. We use Chebfun's \texttt{aaa} function \autocite{nakatsukasa2018aaa, driscoll2014chebfun}. For each time $t_j$, an AAA approximant $r_j(x)$ is constructed from the data $\{(x_i, u_{ij})\}_{i=1}^N$. The \texttt{aaa} function returns, in addition to the approximant $r_j(x)$, a vector of poles $\vb p_j$ and corresponding residues $\vb r_j$. We remove the poles whose corresponding residue is less than $10^{-4}$ in modulus, leaving the modified vectors $\tilde{\vb p}_j$ and $\tilde{\vb r}_j$ for the poles and residues, respectively. More constraints on the poles may be needed for more general nonlinear pdes, such as constraining the poles to be purely imaginary, though we do not investigate these extensions here.

Now that we have the vectors of poles $\{\vb p_j\}_{j=1,\ldots,M}$, we can track the poles of our solution in time.  The procedure works by backtracking in time, starting with an accurate guess $z_0$ for the closest pole to the real line at the time $t_M=2$.  In particular, if we let $\zeta_j$ denote the estimated closest pole to the real line at the time $t_j$, then we set $\zeta_M = z_0$, and we then need to find $\zeta_{M-1},\zeta_{M-2},\ldots,\zeta_1$, in that order. If we have a given $\zeta_j$, then to find $\zeta_{j-1}$ we find the pole in $\tilde{\vb p}_{j-1}$ that is closest to $\zeta_j$, i.e. $\zeta_{j-1} = \operatorname{argmin}_{p_{j-1, \ell}} |p_{j-1,\ell} - \zeta_j|$, where $p_{j-1, \ell}$ is the $\ell$th pole in $\tilde{\vb p}_{j-1}$. This procedure continues up until $j=1$, at which point we have our pole trajectory $\{(t_j, \zeta_j)\}_{j=1}^M$.
\section{Further details of tracking the closest pole with the method of steepest descents}\label{sec:appendixsteepest}

In this section, we give further details on computing \eqref{eq:saddlepointasymptoticexpansionsmallmu}. Using the tangent line $s = s_j + r \e^{\i\delta_j}$, letting $\delta_j$ be the angle in \eqref{eq:saddlepointangles} such that $\cos(\delta_j)>0$, we can integrate in a small $\varepsilon$-neighbourhood of the tangent line to write the contribution at saddle point $s_j$ as \begin{equation}
D_j \sim \int_{-\varepsilon}^\varepsilon \e^{h\left(s_j + r \e^{\i\delta_j}\right)/\mu} \e^{\i\delta_j} \dr.
\end{equation}
Now expand \begin{equation}
h\left(s_j + r\e^{\i\delta_j}\right) \sim h(s_j) + \frac{h''(s_j)\e^{2\i\delta_j}}{2}r^2,\quad r \to 0.
\end{equation}
Thus, extending the integral limits to $\mathbb R$,
\begin{equation}
D_j \sim \e^{h(s_j)/\mu}\e^{\i\delta_j}\int_{-\infty}^\infty \e^{h''(s_j)\e^{2\i\delta_j}r^2/(2\mu)} \dr.
\end{equation}
Now notice that $h''(s_j)\e^{2\i\delta_j} = |h''(s_j)|\e^{\i(\alpha_j + 2\delta_j)}$. We recall that from \eqref{eq:saddlethetacondition} that $\sin(2\delta_j+\alpha_j) = 0$ and $\cos(2\delta_j+\alpha_j)< 0$, and thus $\cos(2\delta_j + \alpha_j) = -1$ so that $\e^{\i(\alpha+j+2\delta_j)} = -1$, giving $h''(s_j)\e^{2\i\delta_j} = -|h''(s_j)|$. Hence,
\begin{equation}\label{eq:individualcomplexsaddlecontribution}
D_j \sim \e^{h(s_j)/\mu}\e^{\i\delta_j}\int_{-\infty}^\infty \e^{-|h''(s_j)|r^2/(2\mu)} \dr = \e^{h(s_j)/\mu}\e^{\i\delta_j}\sqrt{\frac{2\mu\pi}{|h''(s_j)|}},\quad \mu \to 0^+.
\end{equation}
We could sum this over both saddle points, say $s_{j_1}$ and $s_{j_2}$, which contribute to the asymptotic expansion of $D$, leaving $D(\mu) \sim D_{j_1} + D_{j_2}$ as $\mu \to 0^+$.

\begin{figure}[h!]
\centering
\includegraphics[width=0.8\textwidth]{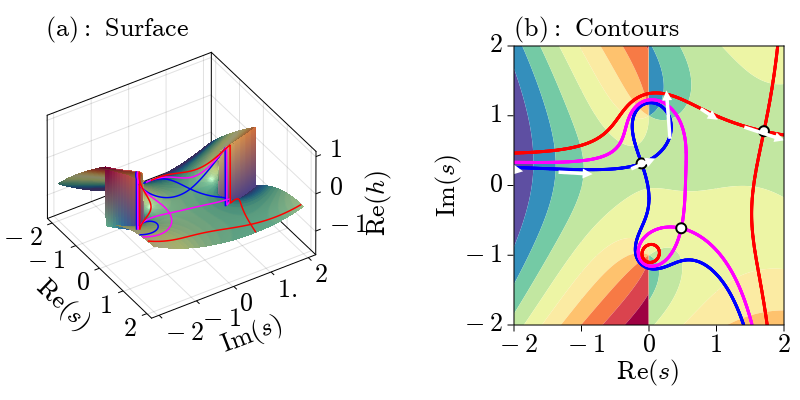}
\caption{Contour plot and surface for the saddle point problem at $z = 2.0713 + 0.48208\i$ and $t = 2$. The colour shows $\Re h(s)$, with $h(s)$ defined in \eqref{eq:denominatordmudependence}, the black-white markers show the locations of the saddle points, and the coloured contours are the contour lines $\Im h(s) = \text{constant}$, where the constant is determined by the imaginary part of $h$ at the saddle point. The surface in (a) is $\Re h(s)$. The arrows in (b) show the path of integration used to go through each saddle point.}\label{fig:saddlepointcontoursurface}
\end{figure}

We now consider the geometry of these saddle-point contours and verify that we can always deform onto them without issue so that \eqref{eq:individualcomplexsaddlecontribution} is valid. For exposition's sake we consider a specific example with $z = 2.0713 + 0.48208\i$, $t = 2$, and $\mu = 0.05$. This point $z$ is approximately the position of the closest pole of $u$ to the real line at $t=2$ for this value of $\mu$. In figure \ref{fig:saddlepointcontoursurface} we show the contour plot and surface for $h(s)$, labelling the saddle points $(s_1, s_2, s_3) \approx (1.703+0.7781\i, -0.10944+0.31644\i, 0.4778-0.61306\i)$.

The contours in figure \ref{fig:saddlepointcontoursurface}(b) appear to show a feature that should not be possible for a saddle-point contour: they loop around. By looking at the surface in (a), we see that the lines appear to rise up a large `cliff' and come back around. In reality, these ``loops'' are simply contours that are going onto another sheet of the Riemann surface defined by $\arctan(z)$, recalling that $\arctan(z)$ has the branch cut $(\e^{-\i\pi/2}\infty, -\i] \cup [\i, \e^{\i\pi/2}\infty)$.   These loops are typical in all plots like figure~\ref{fig:saddlepointcontoursurface}(b), at least for those $z$ that are close to poles of $u$.  Ultimately, we avoid dealing with these loops and the branch cut in particular.  More detail on the structure of this surface and how to deal with the branch cuts is given in the analysis in Chapter 4 of  \textcite{langman2005stokes}.

The colours in the contour plot in figure \ref{fig:saddlepointcontoursurface} can be used to determine which contours lead to divergent integrals. The magenta contour leads to a path with increasing $\Re(h)$, while the blue and red contours lead to paths with decreasing $\Re(h)$. If we let $s_1$ be the saddle point on the red contour, $s_2$ the saddle point on the blue contour, and $s_3$ the saddle point on the magenta contour, then these observations suggest that we want to consider the saddle-point contributions from $s_1$ and $s_2$. For our $z$ we can compute $h(s_1) = -0.5573 - 0.0628\i$, $h(s_2) = -0.5306 - 0.2518\i$, and $h(s_3) = -0.462 - 0.1846\i$. We see that $\Re h(s_1)$ and $\Re h(s_2)$ are indeed close, consistent with the fact that the two saddle-point contributions must have near-equal real parts. Now actually to deform onto these contours, as shown in figure \ref{fig:saddlepointcontoursurface}(b), we first bridge from the real line onto the blue contour, through the saddle point, where we then bridge up to the red contour, through the saddle point, and then eventually back onto the real line and through to infinity. The contributions over the small bridges will be negligible in the limit $\mu\rightarrow 0^+$.

\section{Initial conditions that do not have simple poles} \label{eq:gammanot1}

{
In this section, we outline the key ideas for the small-time analysis for the more general initial condition (\ref{eq:ic2}), where $\beta>0$.  The motivation is to explore the added complications that arise when the initial condition has singularities that are not simple poles, remembering that solutions Burgers' equation must have simple poles with local behaviour given by (\ref{eq:localpole}).

For what follows it is sufficient to suppose that
\begin{equation}
u_0\sim \frac{A}{(x-\mathrm{i})^\beta}
\quad\mbox{as}\quad x\rightarrow\mathrm{i}.
\label{eq:u0appen}
\end{equation}
In order to match with (\ref{eq:ic2}), we would require $A=(-\mathrm{i}/2)^\beta$; however, the precise value of $A$ is unimportant in the following calculations.  By writing out the naive expansion (\ref{eq:naiveexpansion}), we find
\begin{equation}
u_1\sim \frac{\mu\beta(\beta+1)A}{(x-\mathrm{i})^{\beta+2}}
+\frac{\beta A^2}{(x-\mathrm{i})^{2\beta+1}}
\quad\mbox{as}\quad x\rightarrow\mathrm{i}.
\label{eq:u1appen}
\end{equation}
Crucially, the first term on the right-hand side of (\ref{eq:u1appen}) dominates if $\beta<1$, while the second term dominates if $\beta>1$.  We summarise each case below.  The less complicated borderline case $\beta=1$ was treated in section~\ref{sec:smalltime}.

}


\subsection{$0<\beta<1$, diffusion dominated}

For $\beta<1$, the diffusive term $\mu u''_0$ dominates the advective term $u_0u'_0$ near $x=\mathrm{i}$.  For a fixed time $t\ll 1$, we use (\ref{eq:u0appen}) and (\ref{eq:u1appen}) to see that $u_0$ balances $tu_1$ when
$$
\xi=\frac{x-\mathrm{i}}{t^{1/2}}=\mathcal{O}(1).
$$
Writing (\ref{eq:u0appen}) using this variable $\xi$, we see the inner scaling should be
$$
u=\frac{1}{t^{\beta/2}}\Phi(\xi,t),
$$
so Burgers' equation becomes
\begin{equation}
t\frac{\partial \Phi}{\partial t}-\frac{\beta}{2}\Phi-\frac{1}{2}\xi\frac{\partial \Phi}{\partial \xi}+t^{(1-\beta)/2}\Phi\frac{\partial \Phi}{\partial \xi}
=\mu\frac{\partial^2 \Phi}{\partial\xi^2}.
\label{eq:innerpdegless1}
\end{equation}
Further, by writing $\Phi\sim \Phi_0(\xi)$ as $t\rightarrow 0^+$, we find the leading order inner problem is
\begin{equation}
-\frac{\beta}{2}\Phi_0-\frac{1}{2}\xi\Phi_0' = \mu\Phi''_0,
\label{eq:innerodegless1}
\end{equation}
\begin{equation}
\Phi_0\sim \frac{A}{\xi^\beta}
\quad\mbox{as}\quad \xi\rightarrow -\mathrm{i}\infty.
\label{eq:farfieldgless1}
\end{equation}
Notably, (\ref{eq:innerodegless1}) is linear, as the nonlinear advective term is neglected at this stage (cf. (\ref{eq:innerpdeconstantterm}), which applies for $\beta=1$).

A Liouville-Green (WKB) argument applied to (\ref{eq:innerodegless1}) shows that
$\Phi_0\sim K\xi^{\beta-1}\mathrm{e}^{-\xi^2/4\mu}$ as $\xi\rightarrow \mathrm{i}\infty$.  Alternatively, we can solve the linear problem (\ref{eq:innerodegless1})-(\ref{eq:farfieldgless1}) exactly in terms of parabolic cylinder functions to give
$$
\Phi_0=\frac{A}{(-2\mu)^{\beta/2}}\,\mathrm{e}^{-\xi^2/8\mu}
U(\beta-\sfrac{1}{2},\mathrm{i}\xi/\sqrt{2\mu}),
$$
which provides the same exponential growth up the imaginary axis.  Here $\Phi_0$ is entire and so does not have the singularity structure we are after.  However, as we move up the imaginary $\xi$-axis, the neglected term $t^{(1-\beta)/2}\Phi\Phi_\xi$ in (\ref{eq:innerpdegless1}) eventually becomes the same size as the largest terms in (\ref{eq:innerodegless1}), which means it can no longer be ignored.  Thus the inner expansion breaks down when
\begin{equation}
\xi^{1+\beta}\mathrm{e}^{-\xi^2/4\mu}=
\mathcal{O}(t^{(1-\beta)/2}\xi^{2\beta-1}
\mathrm{e}^{-\xi^2/2\mu})
\quad\mbox{as}\quad \xi\rightarrow\mathrm{i}\infty.
\label{eq:breakdown}
\end{equation}
Solving asymptotically for $\xi$, we arrive at the new, rather complicated, scaling
$$
\xi=\mathrm{i}(2\mu(1-\beta))^{1/2}\ln^{1/2}(1/t)
+\frac{\mathrm{i}(2-\beta)\mu^{1/2}}{(2(1-\beta))^{1/2}}
\frac{\ln\ln(1/t)+\ln(2\mu(1-\beta))}{\ln^{1/2}(1/t)}
+\left(\frac{2\mu}{1-\beta}\right)^{1/2}\frac{X}{\ln^{1/2}(1/t)},
$$
where $X=\mathcal{O}(1)$.

In terms of $X$ and $t$, the new ansatz for the inner variable is
$$
\Phi=\left(\frac{2\mu(1-\beta)\ln(1/t)}{t^{1-\beta}}
\right)^{1/2}F(X,t).
$$
To leading order, we write $F\sim F_0(X)$, so that
\begin{equation}
-\mathrm{i}F_0'+2F_0F'_0=F''_0,
\label{eq:rescaleodegless1}
\end{equation}
subject to a matching condition of the form $F_0\sim \mathrm{i}^{\beta-1}K\mathrm{e}^{-\mathrm{i}X}$ as $X\rightarrow -\mathrm{i}\infty$, where the constant $K=\mathrm{i}A\sqrt{2\pi}(2\mu)^{1/2-\beta}/\Gamma(\beta)$ is related to the far field behaviour of $\Phi_0$ above.  The solution of (\ref{eq:rescaleodegless1}) is
\begin{equation}
F_0=\frac{\mathrm{i}}
{1-\mathrm{e}^{\mathrm{i}X}/\mathrm{i}^\beta K},
\label{eq:solnF0}
\end{equation}
which has simple poles at $X=-\mathrm{i}\log (\mathrm{i}^\beta K)+2\pi n$, where $n$ is an integer.

To illustrate these results further, suppose that $\beta=1/2$.  In this case, the constants $A$ and $K$ are $A=(1-\mathrm{i})/2$ and $K=\mathrm{i}^{1/2}$.  The solution for $F_0$ becomes $F_0=\mathrm{i}/(1+\mathrm{i}\mathrm{e}^{\mathrm{i}X})$, with infinitely many simple poles at $X=\pi/2+2\pi n$. Therefore, in terms of the original variables, the poles emerge from $x=\mathrm{i}$ as
\begin{equation}
s(t)\sim \mathrm{i}+\mathrm{i}\mu^{1/2}t^{1/2}\ln^{1/2}(1/t)
+\frac{3\mathrm{i}\mu^{1/2}t^{1/2}(\ln\ln(1/t)+\ln\mu)}
{2\ln^{1/2}(1/t)}
+\frac{\mu^{1/2}\pi(1+4n)t^{1/2}}{\ln^{1/2}(1/t)}
\quad\mbox{as}\quad t\rightarrow 0^+.
\label{eq:gamma12poles}
\end{equation}
By carefully unpacking the rescalings, we can show that these simple poles are of the form $u\sim -2\mu/(x-s(t))$ as required (see (\ref{eq:localpole})).  Supporting figures for $\beta=1/2$ are provided in the Supplementary Material.

Unfortunately, the asymptotic form (\ref{eq:gamma12poles}) for the location of the poles itself breaks down when $|n|=\mathcal{O}(\ln(1/t))$ for integer $n$,  which requires a further treatment that we shall not pursue here, other than the following remark: the above analysis implies (see (\ref{eq:breakdown})) non-uniformity at
$$
\mathrm{e}^{\xi^2/4\mu}=\mathcal{O}\left(
t^{(1-\beta)/2}
\right),
$$
ignoring the pre-exponential factors that can be neglected for the purposes of the current calculation -- hence
\begin{equation}
\xi^2\sim -2\mu(1-\beta)\ln(1/t)+2n\pi\mathrm{i};
\label{eq:simplifiedscaling}
\end{equation}
the above discussion details the case $n=\mathcal{O}(1)$ as $t\rightarrow 0^+$, leading to (\ref{eq:solnF0}), but (\ref{eq:simplifiedscaling}) confirms how the poles again approach alignment $\mathrm{arg}(\xi)=\pi/4$, $3\pi/4$, for $|n|=\mathcal{O}(\ln(1/t))$.

\subsection{$\beta>1$, advection dominated}

In the other case, $\beta>1$, the advective term $u_0u'_0$ dominates the diffusive term $\mu u''_0$ near $x=\mathrm{i}$.  Here, for a fixed small time, we see that $u_0$ balances the next term in (\ref{eq:naiveexpansion}) when
$$
\xi=\frac{x-\mathrm{i}}{t^{1/(\beta+1)}}=\mathcal{O}(1);
$$
thus, the required inner scaling for $u$ is
$$
u=\frac{1}{t^{\beta/(\beta+1)}}\Phi(\xi,t).
$$
Burgers' equation in these variables becomes
\begin{equation}
t\frac{\partial \Phi}{\partial t}-\frac{\beta}{\beta+1}\Phi
-\frac{1}{\beta+1}\xi\frac{\partial \Phi}{\partial \xi}
+\Phi\frac{\partial \Phi}{\partial \xi}
=\mu t^{(\beta-1)/(\beta+1)}\frac{\partial^2 \Phi}{\partial\xi^2},
\label{eq:innerpdeggreater1}
\end{equation}
which, to leading order with $\Phi\sim\Phi_0(\xi)$, becomes
\begin{equation}
-\frac{\beta}{\beta+1}\Phi_0-\frac{1}{\beta+1}\xi\Phi_0'
+\Phi_0\Phi'_0 = 0,
\label{eq:innerodeggreater1}
\end{equation}
together with (\ref{eq:farfieldgless1}).  Here the inner problem (\ref{eq:innerodeggreater1}) is first order and nonlinear, as the diffusive term does not appear.

We may solve (\ref{eq:innerodeggreater1}) with (\ref{eq:farfieldgless1}) exactly to give
$$
\Phi_0(\xi-\Phi_0)^\beta=A.
$$
Analysis of this implicit solution shows that $\Phi_0$ has square-root branch points at $\xi=\xi_0$ of the form
\begin{equation}
\Phi_0\sim \frac{\xi_0}{\beta+1}
-\frac{(2\beta)^{1/2}\xi_0^{1/2}}{\beta+1}
(\xi-\xi_0)^{1/2}
\quad\mbox{as}\quad \xi\rightarrow\xi_0,
\label{eq:phi0ggreater1}
\end{equation}
where the location of the singularities are related to $A$ and $\beta$ via $\xi_0=(1+\beta)(A/\beta^\beta)^{1/(\beta+1)}$.

Now viscous Burgers' equation does not have branch-point singularities, which suggests that this analysis must break down where the as-yet-neglected diffusion term kicks in.  By substituting (\ref{eq:phi0ggreater1}) back into (\ref{eq:innerpdeggreater1}), we see that the diffusion terms is now equally important when
$$
X=\frac{\xi-\xi_0}{t^{2(\beta-1)/3(\beta+1)}},
$$
which leads to the new scaling
$$
\Phi=\frac{\xi_0}{\beta+1}
-\frac{(2\beta)^{1/2}\xi_0^{1/2}}{\beta+1}
t^{(\beta-1)/3(\beta+1)}F(X,t).
$$
By writing $F\sim F_0(X)+t^{(\beta-1)/3(\beta+1)}F_1(X)$, we find that, to leading order,
$$
\frac{(\beta\xi_0)^{1/2}}{\beta+1}(-1+2F_0F'_0)=-\sqrt{2}\mu F''_0.
$$
This second-order equation is nonlinear, but we can integrate directly to give a Riccati equation which can be solved exactly. The solution in terms of Airy functions that satisfies the far-field condition is
$$
F_0=-\frac{\mathrm{Ai}'(\lambda^{2/3}X)}
{\lambda^{1/3}\mathrm{Ai}(\lambda^{2/3}X)},
$$
where $\lambda=-(\beta\xi_0)^{1/2}/\sqrt{2}\mu(\beta+1)$.  The poles for $F_0$ are zeros of $\mathrm{Ai}(\lambda^{2/3}X)$.

To take an example, say $\beta=2$.  In that case, the implicit solution for $\Phi_0$ is $\Phi_0(\xi-\Phi_0)^2=-1/4$.  Only one of the three solution has the appropriate far-field behaviour $\Phi_0\sim -1/4\xi^2$ (the other two have $\Phi_0\sim \xi$).  The three branch points are at $\xi_0=-3/16^{1/3}$, $3\mathrm{e}^{\mathrm{i}\pi/3}/16^{1/3}$, $3\mathrm{e}^{-\mathrm{i}\pi/3}/16^{1/3}$; the branch cuts are arbitrary at this stage.  This inner scaling breaks down near each of the three branch points, where a new inner variable is $X=(\xi-\xi_0)/t^{2/9}$.  If we denote the zeros of the Airy function by $\lambda^{2/3}X_0\in \mathbb R$, in terms of the original variable, the poles emerge from $x=\mathrm{i}$ like
$$
s(t)\sim \mathrm{i}+t^{1/3}\xi_0+t^{5/9}X_0.
$$
For $\beta=2$, we have $\lambda^{2/3}=(-1/16)^{1/9}/(3\mu^2)^{1/3}$, so care needs to be taken with the appropriate choice of $\lambda$, which will also dictate the direction the branch cuts of $\Phi_0$ take from each of the $\xi_0$.  The precise form of the branch cuts could be determined by applying a type of Rankine–Hugoniot shock condition in the complex plane, but we shall not pursue these ideas here.  Regardless, this asymptotic approximation will break down as $|X_0|$ increases, requiring an even more refined treatment, which we shall not describe.  We have included numerical solutions for $\beta=2$ in the Supplementary Material.



\end{document}